\documentclass[11pt,a4paper]{article}

\usepackage{epsf,epsfig,amsfonts,amsgen,amsmath,amstext,amsbsy,amsopn,amsthm
}
\usepackage{amsmath,times,mathptmx}
\usepackage{amsfonts,amsthm,amssymb}
\usepackage{amsfonts}
\usepackage{graphics}
\usepackage{latexsym,bm}
\usepackage{amsfonts,amsthm,amssymb,bbding}
\usepackage{indentfirst}
\usepackage{graphicx}
\usepackage{color}
\usepackage[colorlinks=true,anchorcolor=blue,filecolor=blue,linkcolor=blue,urlcolor=blue,citecolor=blue]{hyperref}
\usepackage{float}
\usepackage{tikz}
\newtheorem{thm}{Theorem}

\newtheorem{lemma}{Lemma}
\newtheorem{false statement}{False statement}

\theoremstyle{definition}

\newtheorem{claim}{Claim}

\newtheorem{remark}[claim]{Remark}

\begin{document}

\title{Maximizing the index of signed complete graphs with spanning trees on $k$ pendant vertices
\footnote{Supported by NSFC (Nos. 12361071 and 11901498).}}
\author{Dan Li, Minghui Yan, Zhaolin Teng\thanks{Corresponding author. E-mail: zlteng03@163.com.}\\
{\footnotesize College of Mathematics and System Science, Xinjiang University, Urumqi 830046, PR China}
}
\date{}

\maketitle {\flushleft\large\bf Abstract:}
A signed graph $\Sigma=(G,\sigma)$ consists of an underlying graph $G=(V,E)$ with a sign function $\sigma:E\rightarrow\{-1,1\}$. 
Let $A(\Sigma)$ be the adjacency matrix of $\Sigma$ and $\lambda_1(\Sigma)$ denote the largest eigenvalue (index) of $\Sigma$.
Define $(K_n,H^-)$ as a signed complete graph whose negative edges induce a subgraph $H$. In this paper, we focus on the following problem:
which spanning tree $T$ with a given number of pendant vertices makes the $\lambda_1(A(\Sigma))$ of the unbalanced $(K_n,T^-)$ as large as possible?
To answer the problem, we characterize the extremal signed graph with maximum $\lambda_1(A(\Sigma))$ among graphs of type $(K_n,T^-)$.

\vspace{0.1cm}
\begin{flushleft}
\end{flushleft}

\section{Introduction}
Throughout this paper, unless stated otherwise, all the underlying graphs are considered as simple and connected.
A signed graph $\Sigma=(G,\sigma)$ consists of an underlying graph $G=(V,E)$ with a sign function $\sigma:E\rightarrow\{-1,1\}$.
The (unsigned) graph $G$ is said to be the underlying graph of $\Sigma$, and the function $\sigma$ is called the signature of $\Sigma$.
The degree of a vertex $v$, denoted by $d(v)$, is the number of the edges incident with $v$. Especially, a vertex of degree one is called pendant vertex.
In signed graphs, edge signs are usually interpreted as $\pm1$. An edge $e$ is positive (resp. negative) if $\sigma(e)=+1$ (resp. $\sigma(e)=-1$). 
The sign of a signed cycle $C$ is given by $\sigma(C)=\prod_{e\in C}\sigma(e)$. A cycle $C$ is positive (resp. negative) if its sign is $+$ (resp. $-$); alternatively, we can say that a signed cycle is positive (resp. negative) if it contains an even (resp. odd) number of negatives edges.
A signed graph is balanced if there are no negative cycles; otherwise it is unbalanced, which
was first introduced in works of Harary \cite{F. Harary} and Cartwright and Harary \cite{Cartwright D}.
And the matroids of graphs were extended to matroids of signed graphs
by Zaslavsky \cite{T. Zaslavsky}. Chaiken \cite{S. Chaiken} and Zaslavsky \cite{T. Zaslavsky} obtained the Matrix-Tree Theorem for signed graph independently. In fact, the theory of signed graphs is a special case of that of gain graphs and of biased graphs \cite{T. Zaslavsky1}.
At the very beginning, signed graphs were studied in the context of social psychology, where the vertices are considered as individuals, and the positive edges and the negative edges represent friendships and enmities between them, respectively.
Note that unsigned graphs are treated as (balanced) signed graphs whose edges all get positive signs, that is, the all-positive signature.

The adjacency matrix of $\Sigma$ is defined as $A(\Sigma)=(a^{\sigma}_{ij})$, where $a^{\sigma}_{ij}=\sigma(v_iv_j)$ if $v_i \thicksim v_j$,  and $0$ otherwise. 
Let $M=M(\Sigma)$ be a real symmetric graph matrix of a signed graph $\Sigma=(G,\sigma)$ and $P_{M}(\lambda)=\mbox{det}(\lambda I-M)$ be the $M$-characteristic polynomial, where $I$ is the identity matrix. The spectrum of $M$ is called the $M$-spectrum of the signed graph $\Sigma$. 
Note that $A(\Sigma)$ is a real symmetric matrix.
As usual, we use
$\lambda_1(\Sigma) \geq\lambda_2(\Sigma) \geq\cdots  \geq\lambda_n(\Sigma)$
to denote the spectrum of $\Sigma$, where $\lambda_n(\Sigma)$ is the least eigenvalue and $\lambda_1(\Sigma)$ is the largest eigenvalue that is often called the index of $\Sigma$.
Then the spectral radius of $\Sigma$ is the largest absolute value of the eigenvalues of $A(\Sigma)$, denoted by $\rho(A(\Sigma))$, i.e., $\rho(A(\Sigma))=\max\{\lambda_1(A(\Sigma)),-\lambda_n(A(\Sigma))\}$.

Extremal graph theory deals with the problem of determining extremal values and extremal graphs for a given graph invariant in a given set of graphs.
In \cite{B.F. Wu}, Wu, Xiao and Hong determined the unique tree of order $n$ with $k$ pendant vertices which achieved the maximal spectral radius. They showed that the maximal spectral radius was obtained uniquely at $T_{n,k}$, where $T_{n,k}$ was a tree of order $n$ obtained from a star $K_{1,k}$ and $k$ paths of almost equal lengths by joining each pendant vertex of $K_{1,k}$ to an end vertex of one path. Guo \cite{S.G. Guo} determined the graphs with the largest spectral radius among all the unicyclic and bicyclic graphs with $n$ vertices and $k$ pendant vertices, respectively.
Petrovi\'c and Borovi\'canin \cite{M. Petrovic} then determined the unique graph with largest spectral radius among all tricyclic graphs with $n$ vertices and $k$ pendant vertices. 
A graph $G$ is called a cactus if any two of its cycles have at most one common vertex. In \cite{J.C. Wu}, Wu, Deng and Jiang identified the graphs with the largest spectral radius among all the cacti with $n$ vertices and $k$ pendant vertices.
Apart from these, there are also some research results on the least eigenvalues of graphs with given pendant vertices, which can be seen in \cite{R.F. Liu, R.F. Liu1, R.D. Xing}.
In addition, relations between the index or the least eigenvalue and a given graph invariant in a given set of signed graphs have been studied. Some recent research can be found in \cite{M. Brunetti, M. Souri, Z.L. Teng, M.H. Yan}.

In this paper, we deal with signed graphs with extremal eigenvalues.
It is worth noting that Koledin and Stani\'{c} \cite{T. Koledin} studied the connected signed graphs of fixed order, size, and the number of negative edges with maximum index. They conjectured that if $\Sigma$ is a signed complete graph of order $n$ with $k$ negative edges, $k<n-1$ and $\Sigma$ has maximum index, then the negative edges induce the signed star $K_{1,k}$,
which was proved by Akbari, Dalvandi, Heydari and Maghasedi \cite{S. Akbari} for signed complete graphs whose negative edges form a tree.
And the conjecture in \cite{T. Koledin} was then completely confirmed by Ghorbani and Majidi \cite{E. Ghorbani}.
Motivated by these, we have characterized the extremal signed graph with maximum $\lambda_1(A(\Sigma))$ and minimum $\lambda_n(A(\Sigma))$ among graphs of type $(K_n,T^-)$, respectively in \cite{D. Li}.
In this article, we continue the line of work in \cite{E. Ghorbani, T. Koledin} on identifying the complete signed graphs with maximal index comprising $n$ vertices and $m$ negative edges.
There have also been some good results on this issue.
Kafai, Heydari, Rad and Maghasedi \cite{N. Kafai} showed that among all signed complete graphs of order $n>5$ whose negative edges induce a unicyclic graph of order $k$ and maximizes the index, the negative edges induce a triangle with all remaining vertices being pendant at the same vertex of the triangle.
In \cite{E. Ghorbani1}, Ghorbani and Majidi identified the unique signed graph with maximal index among complete signed graphs whose negative edges induced a tree of diameter at least $d$ for any given $d$. They also gave the smallest minimum eigenvalue of complete signed graphs with $n$ vertices whose negative edges induced a tree.
We then have the idea that take some other given parameters into consideration. 

Let $K_n$ and $(K_n,H^-)$ denote the complete graph of order $n$ and a signed complete graph whose negative edges induce a subgraph $H$, respectively. Let $T$ be a spanning tree of $K_n$.
In this paper, we mainly focus on the problem: which spanning tree $T$ with a given number of pendant vertices makes the $\lambda_1(A(\Sigma))$ of the unbalanced $(K_n,T^-)$ as large as possible? 
Then we have the following result.

\begin{thm} \label{th1} 
Let $T$ be a spanning tree with $k$ pendant vertices of $K_n$ and $n\geq6$. If $\Sigma=(K_n,T^-)$ is an unbalanced signed complete graph with the maximum index, then $T\cong T^*$.
\end{thm}

\begin{center}
\tikzset{every picture/.style={line width=0.75pt}} 

\begin{tikzpicture}[x=0.75pt,y=0.75pt,yscale=-1,xscale=1]

\draw    (134.33,85) -- (174.01,85.21) ;
\draw [shift={(174.01,85.21)}, rotate = 0.31] [color={rgb, 255:red, 0; green, 0; blue, 0 }  ][fill={rgb, 255:red, 0; green, 0; blue, 0 }  ][line width=0.75]      (0, 0) circle [x radius= 1.34, y radius= 1.34]   ;
\draw [shift={(134.33,85)}, rotate = 0.31] [color={rgb, 255:red, 0; green, 0; blue, 0 }  ][fill={rgb, 255:red, 0; green, 0; blue, 0 }  ][line width=0.75]      (0, 0) circle [x radius= 1.34, y radius= 1.34]   ;
\draw  [dash pattern={on 0.84pt off 2.51pt}]  (174.01,85.21) -- (203.01,85.21) ;
\draw [shift={(203.01,85.21)}, rotate = 0] [color={rgb, 255:red, 0; green, 0; blue, 0 }  ][fill={rgb, 255:red, 0; green, 0; blue, 0 }  ][line width=0.75]      (0, 0) circle [x radius= 1.34, y radius= 1.34]   ;
\draw [shift={(174.01,85.21)}, rotate = 0] [color={rgb, 255:red, 0; green, 0; blue, 0 }  ][fill={rgb, 255:red, 0; green, 0; blue, 0 }  ][line width=0.75]      (0, 0) circle [x radius= 1.34, y radius= 1.34]   ;
\draw    (203.01,85.21) -- (240.01,85.21) ;
\draw [shift={(240.01,85.21)}, rotate = 0] [color={rgb, 255:red, 0; green, 0; blue, 0 }  ][fill={rgb, 255:red, 0; green, 0; blue, 0 }  ][line width=0.75]      (0, 0) circle [x radius= 1.34, y radius= 1.34]   ;
\draw [shift={(203.01,85.21)}, rotate = 0] [color={rgb, 255:red, 0; green, 0; blue, 0 }  ][fill={rgb, 255:red, 0; green, 0; blue, 0 }  ][line width=0.75]      (0, 0) circle [x radius= 1.34, y radius= 1.34]   ;
\draw    (240.01,85.21) -- (221.01,59.21) ;
\draw [shift={(221.01,59.21)}, rotate = 233.84] [color={rgb, 255:red, 0; green, 0; blue, 0 }  ][fill={rgb, 255:red, 0; green, 0; blue, 0 }  ][line width=0.75]      (0, 0) circle [x radius= 1.34, y radius= 1.34]   ;
\draw [shift={(240.01,85.21)}, rotate = 233.84] [color={rgb, 255:red, 0; green, 0; blue, 0 }  ][fill={rgb, 255:red, 0; green, 0; blue, 0 }  ][line width=0.75]      (0, 0) circle [x radius= 1.34, y radius= 1.34]   ;
\draw    (236.01,59.21) -- (240.01,85.21) ;
\draw [shift={(240.01,85.21)}, rotate = 81.25] [color={rgb, 255:red, 0; green, 0; blue, 0 }  ][fill={rgb, 255:red, 0; green, 0; blue, 0 }  ][line width=0.75]      (0, 0) circle [x radius= 1.34, y radius= 1.34]   ;
\draw [shift={(236.01,59.21)}, rotate = 81.25] [color={rgb, 255:red, 0; green, 0; blue, 0 }  ][fill={rgb, 255:red, 0; green, 0; blue, 0 }  ][line width=0.75]      (0, 0) circle [x radius= 1.34, y radius= 1.34]   ;
\draw    (258.01,59.21) -- (240.01,85.21) ;
\draw [shift={(240.01,85.21)}, rotate = 124.7] [color={rgb, 255:red, 0; green, 0; blue, 0 }  ][fill={rgb, 255:red, 0; green, 0; blue, 0 }  ][line width=0.75]      (0, 0) circle [x radius= 1.34, y radius= 1.34]   ;
\draw [shift={(258.01,59.21)}, rotate = 124.7] [color={rgb, 255:red, 0; green, 0; blue, 0 }  ][fill={rgb, 255:red, 0; green, 0; blue, 0 }  ][line width=0.75]      (0, 0) circle [x radius= 1.34, y radius= 1.34]   ;
\draw  [dash pattern={on 0.84pt off 2.51pt}]  (236.01,59.21) -- (258.01,59.21) ;
\draw [shift={(258.01,59.21)}, rotate = 0] [color={rgb, 255:red, 0; green, 0; blue, 0 }  ][fill={rgb, 255:red, 0; green, 0; blue, 0 }  ][line width=0.75]      (0, 0) circle [x radius= 1.34, y radius= 1.34]   ;
\draw [shift={(236.01,59.21)}, rotate = 0] [color={rgb, 255:red, 0; green, 0; blue, 0 }  ][fill={rgb, 255:red, 0; green, 0; blue, 0 }  ][line width=0.75]      (0, 0) circle [x radius= 1.34, y radius= 1.34]   ;
\draw   (134.83,88) .. controls (134.83,92.67) and (137.16,95) .. (141.83,95) -- (176.83,95) .. controls (183.5,95) and (186.83,97.33) .. (186.83,102) .. controls (186.83,97.33) and (190.16,95) .. (196.83,95)(193.83,95) -- (231.83,95) .. controls (236.5,95) and (238.83,92.67) .. (238.83,88) ;
\draw   (257.83,56) .. controls (257.83,51.33) and (255.5,49) .. (250.83,49) -- (249.32,49) .. controls (242.65,49) and (239.32,46.67) .. (239.32,42) .. controls (239.32,46.67) and (235.99,49) .. (229.32,49)(232.32,49) -- (227.83,49) .. controls (223.16,49) and (220.83,51.33) .. (220.83,56) ;
\draw    (437.01,60.21) -- (418.01,34.21) ;
\draw [shift={(418.01,34.21)}, rotate = 233.84] [color={rgb, 255:red, 0; green, 0; blue, 0 }  ][fill={rgb, 255:red, 0; green, 0; blue, 0 }  ][line width=0.75]      (0, 0) circle [x radius= 1.34, y radius= 1.34]   ;
\draw [shift={(437.01,60.21)}, rotate = 233.84] [color={rgb, 255:red, 0; green, 0; blue, 0 }  ][fill={rgb, 255:red, 0; green, 0; blue, 0 }  ][line width=0.75]      (0, 0) circle [x radius= 1.34, y radius= 1.34]   ;
\draw    (433.01,34.21) -- (437.01,60.21) ;
\draw [shift={(437.01,60.21)}, rotate = 81.25] [color={rgb, 255:red, 0; green, 0; blue, 0 }  ][fill={rgb, 255:red, 0; green, 0; blue, 0 }  ][line width=0.75]      (0, 0) circle [x radius= 1.34, y radius= 1.34]   ;
\draw [shift={(433.01,34.21)}, rotate = 81.25] [color={rgb, 255:red, 0; green, 0; blue, 0 }  ][fill={rgb, 255:red, 0; green, 0; blue, 0 }  ][line width=0.75]      (0, 0) circle [x radius= 1.34, y radius= 1.34]   ;
\draw    (455.01,34.21) -- (437.01,60.21) ;
\draw [shift={(437.01,60.21)}, rotate = 124.7] [color={rgb, 255:red, 0; green, 0; blue, 0 }  ][fill={rgb, 255:red, 0; green, 0; blue, 0 }  ][line width=0.75]      (0, 0) circle [x radius= 1.34, y radius= 1.34]   ;
\draw [shift={(455.01,34.21)}, rotate = 124.7] [color={rgb, 255:red, 0; green, 0; blue, 0 }  ][fill={rgb, 255:red, 0; green, 0; blue, 0 }  ][line width=0.75]      (0, 0) circle [x radius= 1.34, y radius= 1.34]   ;
\draw  [dash pattern={on 0.84pt off 2.51pt}]  (433.01,34.21) -- (455.01,34.21) ;
\draw [shift={(455.01,34.21)}, rotate = 0] [color={rgb, 255:red, 0; green, 0; blue, 0 }  ][fill={rgb, 255:red, 0; green, 0; blue, 0 }  ][line width=0.75]      (0, 0) circle [x radius= 1.34, y radius= 1.34]   ;
\draw [shift={(433.01,34.21)}, rotate = 0] [color={rgb, 255:red, 0; green, 0; blue, 0 }  ][fill={rgb, 255:red, 0; green, 0; blue, 0 }  ][line width=0.75]      (0, 0) circle [x radius= 1.34, y radius= 1.34]   ;
\draw    (437.01,60.21) -- (406.37,77.47) ;
\draw [shift={(406.37,77.47)}, rotate = 150.61] [color={rgb, 255:red, 0; green, 0; blue, 0 }  ][fill={rgb, 255:red, 0; green, 0; blue, 0 }  ][line width=0.75]      (0, 0) circle [x radius= 1.34, y radius= 1.34]   ;
\draw [shift={(437.01,60.21)}, rotate = 150.61] [color={rgb, 255:red, 0; green, 0; blue, 0 }  ][fill={rgb, 255:red, 0; green, 0; blue, 0 }  ][line width=0.75]      (0, 0) circle [x radius= 1.34, y radius= 1.34]   ;
\draw    (437.01,60.21) -- (426,78) ;
\draw [shift={(426,78)}, rotate = 121.77] [color={rgb, 255:red, 0; green, 0; blue, 0 }  ][fill={rgb, 255:red, 0; green, 0; blue, 0 }  ][line width=0.75]      (0, 0) circle [x radius= 1.34, y radius= 1.34]   ;
\draw [shift={(437.01,60.21)}, rotate = 121.77] [color={rgb, 255:red, 0; green, 0; blue, 0 }  ][fill={rgb, 255:red, 0; green, 0; blue, 0 }  ][line width=0.75]      (0, 0) circle [x radius= 1.34, y radius= 1.34]   ;
\draw    (457.37,78.47) -- (437.01,60.21) ;
\draw [shift={(437.01,60.21)}, rotate = 221.9] [color={rgb, 255:red, 0; green, 0; blue, 0 }  ][fill={rgb, 255:red, 0; green, 0; blue, 0 }  ][line width=0.75]      (0, 0) circle [x radius= 1.34, y radius= 1.34]   ;
\draw [shift={(457.37,78.47)}, rotate = 221.9] [color={rgb, 255:red, 0; green, 0; blue, 0 }  ][fill={rgb, 255:red, 0; green, 0; blue, 0 }  ][line width=0.75]      (0, 0) circle [x radius= 1.34, y radius= 1.34]   ;
\draw    (406.37,103.47) -- (406.37,77.47) ;
\draw [shift={(406.37,77.47)}, rotate = 270] [color={rgb, 255:red, 0; green, 0; blue, 0 }  ][fill={rgb, 255:red, 0; green, 0; blue, 0 }  ][line width=0.75]      (0, 0) circle [x radius= 1.34, y radius= 1.34]   ;
\draw [shift={(406.37,103.47)}, rotate = 270] [color={rgb, 255:red, 0; green, 0; blue, 0 }  ][fill={rgb, 255:red, 0; green, 0; blue, 0 }  ][line width=0.75]      (0, 0) circle [x radius= 1.34, y radius= 1.34]   ;
\draw    (426,104) -- (426,78) ;
\draw [shift={(426,78)}, rotate = 270] [color={rgb, 255:red, 0; green, 0; blue, 0 }  ][fill={rgb, 255:red, 0; green, 0; blue, 0 }  ][line width=0.75]      (0, 0) circle [x radius= 1.34, y radius= 1.34]   ;
\draw [shift={(426,104)}, rotate = 270] [color={rgb, 255:red, 0; green, 0; blue, 0 }  ][fill={rgb, 255:red, 0; green, 0; blue, 0 }  ][line width=0.75]      (0, 0) circle [x radius= 1.34, y radius= 1.34]   ;
\draw    (457.37,104.47) -- (457.37,78.47) ;
\draw [shift={(457.37,78.47)}, rotate = 270] [color={rgb, 255:red, 0; green, 0; blue, 0 }  ][fill={rgb, 255:red, 0; green, 0; blue, 0 }  ][line width=0.75]      (0, 0) circle [x radius= 1.34, y radius= 1.34]   ;
\draw [shift={(457.37,104.47)}, rotate = 270] [color={rgb, 255:red, 0; green, 0; blue, 0 }  ][fill={rgb, 255:red, 0; green, 0; blue, 0 }  ][line width=0.75]      (0, 0) circle [x radius= 1.34, y radius= 1.34]   ;
\draw  [dash pattern={on 0.84pt off 2.51pt}]  (426,104) -- (457.37,104.47) ;
\draw [shift={(457.37,104.47)}, rotate = 0.86] [color={rgb, 255:red, 0; green, 0; blue, 0 }  ][fill={rgb, 255:red, 0; green, 0; blue, 0 }  ][line width=0.75]      (0, 0) circle [x radius= 1.34, y radius= 1.34]   ;
\draw [shift={(426,104)}, rotate = 0.86] [color={rgb, 255:red, 0; green, 0; blue, 0 }  ][fill={rgb, 255:red, 0; green, 0; blue, 0 }  ][line width=0.75]      (0, 0) circle [x radius= 1.34, y radius= 1.34]   ;
\draw    (163.67,189.67) -- (184.01,189.74) ;
\draw [shift={(184.01,189.74)}, rotate = 0.22] [color={rgb, 255:red, 0; green, 0; blue, 0 }  ][fill={rgb, 255:red, 0; green, 0; blue, 0 }  ][line width=0.75]      (0, 0) circle [x radius= 1.34, y radius= 1.34]   ;
\draw [shift={(163.67,189.67)}, rotate = 0.22] [color={rgb, 255:red, 0; green, 0; blue, 0 }  ][fill={rgb, 255:red, 0; green, 0; blue, 0 }  ][line width=0.75]      (0, 0) circle [x radius= 1.34, y radius= 1.34]   ;
\draw  [dash pattern={on 0.84pt off 2.51pt}]  (184.01,189.74) -- (204.34,189.82) ;
\draw [shift={(204.34,189.82)}, rotate = 0.22] [color={rgb, 255:red, 0; green, 0; blue, 0 }  ][fill={rgb, 255:red, 0; green, 0; blue, 0 }  ][line width=0.75]      (0, 0) circle [x radius= 1.34, y radius= 1.34]   ;
\draw [shift={(184.01,189.74)}, rotate = 0.22] [color={rgb, 255:red, 0; green, 0; blue, 0 }  ][fill={rgb, 255:red, 0; green, 0; blue, 0 }  ][line width=0.75]      (0, 0) circle [x radius= 1.34, y radius= 1.34]   ;
\draw    (204.34,189.82) -- (224.68,189.9) ;
\draw [shift={(224.68,189.9)}, rotate = 0.22] [color={rgb, 255:red, 0; green, 0; blue, 0 }  ][fill={rgb, 255:red, 0; green, 0; blue, 0 }  ][line width=0.75]      (0, 0) circle [x radius= 1.34, y radius= 1.34]   ;
\draw [shift={(204.34,189.82)}, rotate = 0.22] [color={rgb, 255:red, 0; green, 0; blue, 0 }  ][fill={rgb, 255:red, 0; green, 0; blue, 0 }  ][line width=0.75]      (0, 0) circle [x radius= 1.34, y radius= 1.34]   ;
\draw    (224.68,189.9) -- (240.33,184.33) ;
\draw [shift={(240.33,184.33)}, rotate = 340.42] [color={rgb, 255:red, 0; green, 0; blue, 0 }  ][fill={rgb, 255:red, 0; green, 0; blue, 0 }  ][line width=0.75]      (0, 0) circle [x radius= 1.34, y radius= 1.34]   ;
\draw [shift={(224.68,189.9)}, rotate = 340.42] [color={rgb, 255:red, 0; green, 0; blue, 0 }  ][fill={rgb, 255:red, 0; green, 0; blue, 0 }  ][line width=0.75]      (0, 0) circle [x radius= 1.34, y radius= 1.34]   ;
\draw    (224.68,189.9) -- (242.67,189.67) ;
\draw [shift={(242.67,189.67)}, rotate = 359.25] [color={rgb, 255:red, 0; green, 0; blue, 0 }  ][fill={rgb, 255:red, 0; green, 0; blue, 0 }  ][line width=0.75]      (0, 0) circle [x radius= 1.34, y radius= 1.34]   ;
\draw [shift={(224.68,189.9)}, rotate = 359.25] [color={rgb, 255:red, 0; green, 0; blue, 0 }  ][fill={rgb, 255:red, 0; green, 0; blue, 0 }  ][line width=0.75]      (0, 0) circle [x radius= 1.34, y radius= 1.34]   ;
\draw  [dash pattern={on 0.84pt off 2.51pt}]  (242.67,189.67) -- (259.67,189.67) ;
\draw [shift={(259.67,189.67)}, rotate = 0] [color={rgb, 255:red, 0; green, 0; blue, 0 }  ][fill={rgb, 255:red, 0; green, 0; blue, 0 }  ][line width=0.75]      (0, 0) circle [x radius= 1.34, y radius= 1.34]   ;
\draw [shift={(242.67,189.67)}, rotate = 0] [color={rgb, 255:red, 0; green, 0; blue, 0 }  ][fill={rgb, 255:red, 0; green, 0; blue, 0 }  ][line width=0.75]      (0, 0) circle [x radius= 1.34, y radius= 1.34]   ;
\draw    (259.67,189.67) -- (276.67,189.67) ;
\draw [shift={(276.67,189.67)}, rotate = 0] [color={rgb, 255:red, 0; green, 0; blue, 0 }  ][fill={rgb, 255:red, 0; green, 0; blue, 0 }  ][line width=0.75]      (0, 0) circle [x radius= 1.34, y radius= 1.34]   ;
\draw [shift={(259.67,189.67)}, rotate = 0] [color={rgb, 255:red, 0; green, 0; blue, 0 }  ][fill={rgb, 255:red, 0; green, 0; blue, 0 }  ][line width=0.75]      (0, 0) circle [x radius= 1.34, y radius= 1.34]   ;
\draw    (239.67,194.67) -- (224.68,189.9) ;
\draw [shift={(224.68,189.9)}, rotate = 197.64] [color={rgb, 255:red, 0; green, 0; blue, 0 }  ][fill={rgb, 255:red, 0; green, 0; blue, 0 }  ][line width=0.75]      (0, 0) circle [x radius= 1.34, y radius= 1.34]   ;
\draw [shift={(239.67,194.67)}, rotate = 197.64] [color={rgb, 255:red, 0; green, 0; blue, 0 }  ][fill={rgb, 255:red, 0; green, 0; blue, 0 }  ][line width=0.75]      (0, 0) circle [x radius= 1.34, y radius= 1.34]   ;
\draw  [dash pattern={on 0.84pt off 2.51pt}]  (254.65,199.43) -- (239.67,194.67) ;
\draw [shift={(239.67,194.67)}, rotate = 197.64] [color={rgb, 255:red, 0; green, 0; blue, 0 }  ][fill={rgb, 255:red, 0; green, 0; blue, 0 }  ][line width=0.75]      (0, 0) circle [x radius= 1.34, y radius= 1.34]   ;
\draw [shift={(254.65,199.43)}, rotate = 197.64] [color={rgb, 255:red, 0; green, 0; blue, 0 }  ][fill={rgb, 255:red, 0; green, 0; blue, 0 }  ][line width=0.75]      (0, 0) circle [x radius= 1.34, y radius= 1.34]   ;
\draw    (269.63,204.2) -- (254.65,199.43) ;
\draw [shift={(254.65,199.43)}, rotate = 197.64] [color={rgb, 255:red, 0; green, 0; blue, 0 }  ][fill={rgb, 255:red, 0; green, 0; blue, 0 }  ][line width=0.75]      (0, 0) circle [x radius= 1.34, y radius= 1.34]   ;
\draw [shift={(269.63,204.2)}, rotate = 197.64] [color={rgb, 255:red, 0; green, 0; blue, 0 }  ][fill={rgb, 255:red, 0; green, 0; blue, 0 }  ][line width=0.75]      (0, 0) circle [x radius= 1.34, y radius= 1.34]   ;
\draw  [dash pattern={on 0.84pt off 2.51pt}]  (240.33,184.33) -- (255.98,178.77) ;
\draw [shift={(255.98,178.77)}, rotate = 340.42] [color={rgb, 255:red, 0; green, 0; blue, 0 }  ][fill={rgb, 255:red, 0; green, 0; blue, 0 }  ][line width=0.75]      (0, 0) circle [x radius= 1.34, y radius= 1.34]   ;
\draw [shift={(240.33,184.33)}, rotate = 340.42] [color={rgb, 255:red, 0; green, 0; blue, 0 }  ][fill={rgb, 255:red, 0; green, 0; blue, 0 }  ][line width=0.75]      (0, 0) circle [x radius= 1.34, y radius= 1.34]   ;
\draw    (255.98,178.77) -- (271.63,173.2) ;
\draw [shift={(271.63,173.2)}, rotate = 340.42] [color={rgb, 255:red, 0; green, 0; blue, 0 }  ][fill={rgb, 255:red, 0; green, 0; blue, 0 }  ][line width=0.75]      (0, 0) circle [x radius= 1.34, y radius= 1.34]   ;
\draw [shift={(255.98,178.77)}, rotate = 340.42] [color={rgb, 255:red, 0; green, 0; blue, 0 }  ][fill={rgb, 255:red, 0; green, 0; blue, 0 }  ][line width=0.75]      (0, 0) circle [x radius= 1.34, y radius= 1.34]   ;
\draw    (163.67,189.67) -- (148.68,184.9) ;
\draw [shift={(148.68,184.9)}, rotate = 197.64] [color={rgb, 255:red, 0; green, 0; blue, 0 }  ][fill={rgb, 255:red, 0; green, 0; blue, 0 }  ][line width=0.75]      (0, 0) circle [x radius= 1.34, y radius= 1.34]   ;
\draw [shift={(163.67,189.67)}, rotate = 197.64] [color={rgb, 255:red, 0; green, 0; blue, 0 }  ][fill={rgb, 255:red, 0; green, 0; blue, 0 }  ][line width=0.75]      (0, 0) circle [x radius= 1.34, y radius= 1.34]   ;
\draw  [dash pattern={on 0.84pt off 2.51pt}]  (148.68,184.9) -- (133.7,180.14) ;
\draw [shift={(133.7,180.14)}, rotate = 197.64] [color={rgb, 255:red, 0; green, 0; blue, 0 }  ][fill={rgb, 255:red, 0; green, 0; blue, 0 }  ][line width=0.75]      (0, 0) circle [x radius= 1.34, y radius= 1.34]   ;
\draw [shift={(148.68,184.9)}, rotate = 197.64] [color={rgb, 255:red, 0; green, 0; blue, 0 }  ][fill={rgb, 255:red, 0; green, 0; blue, 0 }  ][line width=0.75]      (0, 0) circle [x radius= 1.34, y radius= 1.34]   ;
\draw    (133.7,180.14) -- (118.72,175.37) ;
\draw [shift={(118.72,175.37)}, rotate = 197.64] [color={rgb, 255:red, 0; green, 0; blue, 0 }  ][fill={rgb, 255:red, 0; green, 0; blue, 0 }  ][line width=0.75]      (0, 0) circle [x radius= 1.34, y radius= 1.34]   ;
\draw [shift={(133.7,180.14)}, rotate = 197.64] [color={rgb, 255:red, 0; green, 0; blue, 0 }  ][fill={rgb, 255:red, 0; green, 0; blue, 0 }  ][line width=0.75]      (0, 0) circle [x radius= 1.34, y radius= 1.34]   ;
\draw    (148.02,195.23) -- (163.67,189.67) ;
\draw [shift={(163.67,189.67)}, rotate = 340.42] [color={rgb, 255:red, 0; green, 0; blue, 0 }  ][fill={rgb, 255:red, 0; green, 0; blue, 0 }  ][line width=0.75]      (0, 0) circle [x radius= 1.34, y radius= 1.34]   ;
\draw [shift={(148.02,195.23)}, rotate = 340.42] [color={rgb, 255:red, 0; green, 0; blue, 0 }  ][fill={rgb, 255:red, 0; green, 0; blue, 0 }  ][line width=0.75]      (0, 0) circle [x radius= 1.34, y radius= 1.34]   ;
\draw  [dash pattern={on 0.84pt off 2.51pt}]  (132.37,200.8) -- (148.02,195.23) ;
\draw [shift={(148.02,195.23)}, rotate = 340.42] [color={rgb, 255:red, 0; green, 0; blue, 0 }  ][fill={rgb, 255:red, 0; green, 0; blue, 0 }  ][line width=0.75]      (0, 0) circle [x radius= 1.34, y radius= 1.34]   ;
\draw [shift={(132.37,200.8)}, rotate = 340.42] [color={rgb, 255:red, 0; green, 0; blue, 0 }  ][fill={rgb, 255:red, 0; green, 0; blue, 0 }  ][line width=0.75]      (0, 0) circle [x radius= 1.34, y radius= 1.34]   ;
\draw    (116.72,206.37) -- (132.37,200.8) ;
\draw [shift={(132.37,200.8)}, rotate = 340.42] [color={rgb, 255:red, 0; green, 0; blue, 0 }  ][fill={rgb, 255:red, 0; green, 0; blue, 0 }  ][line width=0.75]      (0, 0) circle [x radius= 1.34, y radius= 1.34]   ;
\draw [shift={(116.72,206.37)}, rotate = 340.42] [color={rgb, 255:red, 0; green, 0; blue, 0 }  ][fill={rgb, 255:red, 0; green, 0; blue, 0 }  ][line width=0.75]      (0, 0) circle [x radius= 1.34, y radius= 1.34]   ;
\draw    (145.68,189.9) -- (163.67,189.67) ;
\draw [shift={(163.67,189.67)}, rotate = 359.25] [color={rgb, 255:red, 0; green, 0; blue, 0 }  ][fill={rgb, 255:red, 0; green, 0; blue, 0 }  ][line width=0.75]      (0, 0) circle [x radius= 1.34, y radius= 1.34]   ;
\draw [shift={(145.68,189.9)}, rotate = 359.25] [color={rgb, 255:red, 0; green, 0; blue, 0 }  ][fill={rgb, 255:red, 0; green, 0; blue, 0 }  ][line width=0.75]      (0, 0) circle [x radius= 1.34, y radius= 1.34]   ;
\draw  [dash pattern={on 0.84pt off 2.51pt}]  (128.68,189.9) -- (145.68,189.9) ;
\draw [shift={(145.68,189.9)}, rotate = 0] [color={rgb, 255:red, 0; green, 0; blue, 0 }  ][fill={rgb, 255:red, 0; green, 0; blue, 0 }  ][line width=0.75]      (0, 0) circle [x radius= 1.34, y radius= 1.34]   ;
\draw [shift={(128.68,189.9)}, rotate = 0] [color={rgb, 255:red, 0; green, 0; blue, 0 }  ][fill={rgb, 255:red, 0; green, 0; blue, 0 }  ][line width=0.75]      (0, 0) circle [x radius= 1.34, y radius= 1.34]   ;
\draw    (111.68,189.9) -- (128.68,189.9) ;
\draw [shift={(128.68,189.9)}, rotate = 0] [color={rgb, 255:red, 0; green, 0; blue, 0 }  ][fill={rgb, 255:red, 0; green, 0; blue, 0 }  ][line width=0.75]      (0, 0) circle [x radius= 1.34, y radius= 1.34]   ;
\draw [shift={(111.68,189.9)}, rotate = 0] [color={rgb, 255:red, 0; green, 0; blue, 0 }  ][fill={rgb, 255:red, 0; green, 0; blue, 0 }  ][line width=0.75]      (0, 0) circle [x radius= 1.34, y radius= 1.34]   ;
\draw    (404.67,189.67) -- (425.01,189.74) ;
\draw [shift={(425.01,189.74)}, rotate = 0.22] [color={rgb, 255:red, 0; green, 0; blue, 0 }  ][fill={rgb, 255:red, 0; green, 0; blue, 0 }  ][line width=0.75]      (0, 0) circle [x radius= 1.34, y radius= 1.34]   ;
\draw [shift={(404.67,189.67)}, rotate = 0.22] [color={rgb, 255:red, 0; green, 0; blue, 0 }  ][fill={rgb, 255:red, 0; green, 0; blue, 0 }  ][line width=0.75]      (0, 0) circle [x radius= 1.34, y radius= 1.34]   ;
\draw  [dash pattern={on 0.84pt off 2.51pt}]  (425.01,189.74) -- (445.34,189.82) ;
\draw [shift={(445.34,189.82)}, rotate = 0.22] [color={rgb, 255:red, 0; green, 0; blue, 0 }  ][fill={rgb, 255:red, 0; green, 0; blue, 0 }  ][line width=0.75]      (0, 0) circle [x radius= 1.34, y radius= 1.34]   ;
\draw [shift={(425.01,189.74)}, rotate = 0.22] [color={rgb, 255:red, 0; green, 0; blue, 0 }  ][fill={rgb, 255:red, 0; green, 0; blue, 0 }  ][line width=0.75]      (0, 0) circle [x radius= 1.34, y radius= 1.34]   ;
\draw    (445.34,189.82) -- (465.68,189.9) ;
\draw [shift={(465.68,189.9)}, rotate = 0.22] [color={rgb, 255:red, 0; green, 0; blue, 0 }  ][fill={rgb, 255:red, 0; green, 0; blue, 0 }  ][line width=0.75]      (0, 0) circle [x radius= 1.34, y radius= 1.34]   ;
\draw [shift={(445.34,189.82)}, rotate = 0.22] [color={rgb, 255:red, 0; green, 0; blue, 0 }  ][fill={rgb, 255:red, 0; green, 0; blue, 0 }  ][line width=0.75]      (0, 0) circle [x radius= 1.34, y radius= 1.34]   ;
\draw    (465.68,189.9) -- (481.33,184.33) ;
\draw [shift={(481.33,184.33)}, rotate = 340.42] [color={rgb, 255:red, 0; green, 0; blue, 0 }  ][fill={rgb, 255:red, 0; green, 0; blue, 0 }  ][line width=0.75]      (0, 0) circle [x radius= 1.34, y radius= 1.34]   ;
\draw [shift={(465.68,189.9)}, rotate = 340.42] [color={rgb, 255:red, 0; green, 0; blue, 0 }  ][fill={rgb, 255:red, 0; green, 0; blue, 0 }  ][line width=0.75]      (0, 0) circle [x radius= 1.34, y radius= 1.34]   ;
\draw    (465.68,189.9) -- (483.67,189.67) ;
\draw [shift={(483.67,189.67)}, rotate = 359.25] [color={rgb, 255:red, 0; green, 0; blue, 0 }  ][fill={rgb, 255:red, 0; green, 0; blue, 0 }  ][line width=0.75]      (0, 0) circle [x radius= 1.34, y radius= 1.34]   ;
\draw [shift={(465.68,189.9)}, rotate = 359.25] [color={rgb, 255:red, 0; green, 0; blue, 0 }  ][fill={rgb, 255:red, 0; green, 0; blue, 0 }  ][line width=0.75]      (0, 0) circle [x radius= 1.34, y radius= 1.34]   ;
\draw  [dash pattern={on 0.84pt off 2.51pt}]  (483.67,189.67) -- (500.67,189.67) ;
\draw [shift={(500.67,189.67)}, rotate = 0] [color={rgb, 255:red, 0; green, 0; blue, 0 }  ][fill={rgb, 255:red, 0; green, 0; blue, 0 }  ][line width=0.75]      (0, 0) circle [x radius= 1.34, y radius= 1.34]   ;
\draw [shift={(483.67,189.67)}, rotate = 0] [color={rgb, 255:red, 0; green, 0; blue, 0 }  ][fill={rgb, 255:red, 0; green, 0; blue, 0 }  ][line width=0.75]      (0, 0) circle [x radius= 1.34, y radius= 1.34]   ;
\draw    (500.67,189.67) -- (517.67,189.67) ;
\draw [shift={(517.67,189.67)}, rotate = 0] [color={rgb, 255:red, 0; green, 0; blue, 0 }  ][fill={rgb, 255:red, 0; green, 0; blue, 0 }  ][line width=0.75]      (0, 0) circle [x radius= 1.34, y radius= 1.34]   ;
\draw [shift={(500.67,189.67)}, rotate = 0] [color={rgb, 255:red, 0; green, 0; blue, 0 }  ][fill={rgb, 255:red, 0; green, 0; blue, 0 }  ][line width=0.75]      (0, 0) circle [x radius= 1.34, y radius= 1.34]   ;
\draw    (480.67,194.67) -- (465.68,189.9) ;
\draw [shift={(465.68,189.9)}, rotate = 197.64] [color={rgb, 255:red, 0; green, 0; blue, 0 }  ][fill={rgb, 255:red, 0; green, 0; blue, 0 }  ][line width=0.75]      (0, 0) circle [x radius= 1.34, y radius= 1.34]   ;
\draw [shift={(480.67,194.67)}, rotate = 197.64] [color={rgb, 255:red, 0; green, 0; blue, 0 }  ][fill={rgb, 255:red, 0; green, 0; blue, 0 }  ][line width=0.75]      (0, 0) circle [x radius= 1.34, y radius= 1.34]   ;
\draw  [dash pattern={on 0.84pt off 2.51pt}]  (495.65,199.43) -- (480.67,194.67) ;
\draw [shift={(480.67,194.67)}, rotate = 197.64] [color={rgb, 255:red, 0; green, 0; blue, 0 }  ][fill={rgb, 255:red, 0; green, 0; blue, 0 }  ][line width=0.75]      (0, 0) circle [x radius= 1.34, y radius= 1.34]   ;
\draw [shift={(495.65,199.43)}, rotate = 197.64] [color={rgb, 255:red, 0; green, 0; blue, 0 }  ][fill={rgb, 255:red, 0; green, 0; blue, 0 }  ][line width=0.75]      (0, 0) circle [x radius= 1.34, y radius= 1.34]   ;
\draw    (510.63,204.2) -- (495.65,199.43) ;
\draw [shift={(495.65,199.43)}, rotate = 197.64] [color={rgb, 255:red, 0; green, 0; blue, 0 }  ][fill={rgb, 255:red, 0; green, 0; blue, 0 }  ][line width=0.75]      (0, 0) circle [x radius= 1.34, y radius= 1.34]   ;
\draw [shift={(510.63,204.2)}, rotate = 197.64] [color={rgb, 255:red, 0; green, 0; blue, 0 }  ][fill={rgb, 255:red, 0; green, 0; blue, 0 }  ][line width=0.75]      (0, 0) circle [x radius= 1.34, y radius= 1.34]   ;
\draw  [dash pattern={on 0.84pt off 2.51pt}]  (481.33,184.33) -- (496.98,178.77) ;
\draw [shift={(496.98,178.77)}, rotate = 340.42] [color={rgb, 255:red, 0; green, 0; blue, 0 }  ][fill={rgb, 255:red, 0; green, 0; blue, 0 }  ][line width=0.75]      (0, 0) circle [x radius= 1.34, y radius= 1.34]   ;
\draw [shift={(481.33,184.33)}, rotate = 340.42] [color={rgb, 255:red, 0; green, 0; blue, 0 }  ][fill={rgb, 255:red, 0; green, 0; blue, 0 }  ][line width=0.75]      (0, 0) circle [x radius= 1.34, y radius= 1.34]   ;
\draw    (496.98,178.77) -- (512.63,173.2) ;
\draw [shift={(512.63,173.2)}, rotate = 340.42] [color={rgb, 255:red, 0; green, 0; blue, 0 }  ][fill={rgb, 255:red, 0; green, 0; blue, 0 }  ][line width=0.75]      (0, 0) circle [x radius= 1.34, y radius= 1.34]   ;
\draw [shift={(496.98,178.77)}, rotate = 340.42] [color={rgb, 255:red, 0; green, 0; blue, 0 }  ][fill={rgb, 255:red, 0; green, 0; blue, 0 }  ][line width=0.75]      (0, 0) circle [x radius= 1.34, y radius= 1.34]   ;
\draw    (476.33,180) -- (465.68,189.9) ;
\draw [shift={(465.68,189.9)}, rotate = 137.09] [color={rgb, 255:red, 0; green, 0; blue, 0 }  ][fill={rgb, 255:red, 0; green, 0; blue, 0 }  ][line width=0.75]      (0, 0) circle [x radius= 1.34, y radius= 1.34]   ;
\draw [shift={(476.33,180)}, rotate = 137.09] [color={rgb, 255:red, 0; green, 0; blue, 0 }  ][fill={rgb, 255:red, 0; green, 0; blue, 0 }  ][line width=0.75]      (0, 0) circle [x radius= 1.34, y radius= 1.34]   ;
\draw    (470.33,178) -- (465.68,189.9) ;
\draw [shift={(465.68,189.9)}, rotate = 111.34] [color={rgb, 255:red, 0; green, 0; blue, 0 }  ][fill={rgb, 255:red, 0; green, 0; blue, 0 }  ][line width=0.75]      (0, 0) circle [x radius= 1.34, y radius= 1.34]   ;
\draw [shift={(470.33,178)}, rotate = 111.34] [color={rgb, 255:red, 0; green, 0; blue, 0 }  ][fill={rgb, 255:red, 0; green, 0; blue, 0 }  ][line width=0.75]      (0, 0) circle [x radius= 1.34, y radius= 1.34]   ;
\draw    (386.68,189.9) -- (404.67,189.67) ;
\draw [shift={(404.67,189.67)}, rotate = 359.25] [color={rgb, 255:red, 0; green, 0; blue, 0 }  ][fill={rgb, 255:red, 0; green, 0; blue, 0 }  ][line width=0.75]      (0, 0) circle [x radius= 1.34, y radius= 1.34]   ;
\draw [shift={(386.68,189.9)}, rotate = 359.25] [color={rgb, 255:red, 0; green, 0; blue, 0 }  ][fill={rgb, 255:red, 0; green, 0; blue, 0 }  ][line width=0.75]      (0, 0) circle [x radius= 1.34, y radius= 1.34]   ;
\draw  [dash pattern={on 0.84pt off 2.51pt}]  (369.68,189.9) -- (386.68,189.9) ;
\draw [shift={(386.68,189.9)}, rotate = 0] [color={rgb, 255:red, 0; green, 0; blue, 0 }  ][fill={rgb, 255:red, 0; green, 0; blue, 0 }  ][line width=0.75]      (0, 0) circle [x radius= 1.34, y radius= 1.34]   ;
\draw [shift={(369.68,189.9)}, rotate = 0] [color={rgb, 255:red, 0; green, 0; blue, 0 }  ][fill={rgb, 255:red, 0; green, 0; blue, 0 }  ][line width=0.75]      (0, 0) circle [x radius= 1.34, y radius= 1.34]   ;
\draw    (352.68,189.9) -- (369.68,189.9) ;
\draw [shift={(369.68,189.9)}, rotate = 0] [color={rgb, 255:red, 0; green, 0; blue, 0 }  ][fill={rgb, 255:red, 0; green, 0; blue, 0 }  ][line width=0.75]      (0, 0) circle [x radius= 1.34, y radius= 1.34]   ;
\draw [shift={(352.68,189.9)}, rotate = 0] [color={rgb, 255:red, 0; green, 0; blue, 0 }  ][fill={rgb, 255:red, 0; green, 0; blue, 0 }  ][line width=0.75]      (0, 0) circle [x radius= 1.34, y radius= 1.34]   ;
\draw  [dash pattern={on 0.84pt off 2.51pt}]  (486.98,170.1) -- (476.33,180) ;
\draw [shift={(476.33,180)}, rotate = 137.09] [color={rgb, 255:red, 0; green, 0; blue, 0 }  ][fill={rgb, 255:red, 0; green, 0; blue, 0 }  ][line width=0.75]      (0, 0) circle [x radius= 1.34, y radius= 1.34]   ;
\draw [shift={(486.98,170.1)}, rotate = 137.09] [color={rgb, 255:red, 0; green, 0; blue, 0 }  ][fill={rgb, 255:red, 0; green, 0; blue, 0 }  ][line width=0.75]      (0, 0) circle [x radius= 1.34, y radius= 1.34]   ;
\draw    (497.63,160.2) -- (486.98,170.1) ;
\draw [shift={(486.98,170.1)}, rotate = 137.09] [color={rgb, 255:red, 0; green, 0; blue, 0 }  ][fill={rgb, 255:red, 0; green, 0; blue, 0 }  ][line width=0.75]      (0, 0) circle [x radius= 1.34, y radius= 1.34]   ;
\draw [shift={(497.63,160.2)}, rotate = 137.09] [color={rgb, 255:red, 0; green, 0; blue, 0 }  ][fill={rgb, 255:red, 0; green, 0; blue, 0 }  ][line width=0.75]      (0, 0) circle [x radius= 1.34, y radius= 1.34]   ;
\draw  [dash pattern={on 0.84pt off 2.51pt}]  (474.98,166.1) -- (470.33,178) ;
\draw [shift={(470.33,178)}, rotate = 111.34] [color={rgb, 255:red, 0; green, 0; blue, 0 }  ][fill={rgb, 255:red, 0; green, 0; blue, 0 }  ][line width=0.75]      (0, 0) circle [x radius= 1.34, y radius= 1.34]   ;
\draw [shift={(474.98,166.1)}, rotate = 111.34] [color={rgb, 255:red, 0; green, 0; blue, 0 }  ][fill={rgb, 255:red, 0; green, 0; blue, 0 }  ][line width=0.75]      (0, 0) circle [x radius= 1.34, y radius= 1.34]   ;
\draw    (479.63,154.2) -- (474.98,166.1) ;
\draw [shift={(474.98,166.1)}, rotate = 111.34] [color={rgb, 255:red, 0; green, 0; blue, 0 }  ][fill={rgb, 255:red, 0; green, 0; blue, 0 }  ][line width=0.75]      (0, 0) circle [x radius= 1.34, y radius= 1.34]   ;
\draw [shift={(479.63,154.2)}, rotate = 111.34] [color={rgb, 255:red, 0; green, 0; blue, 0 }  ][fill={rgb, 255:red, 0; green, 0; blue, 0 }  ][line width=0.75]      (0, 0) circle [x radius= 1.34, y radius= 1.34]   ;
\draw  [dash pattern={on 0.84pt off 2.51pt}]  (271.63,173.2) -- (276.67,189.67) ;
\draw [shift={(276.67,189.67)}, rotate = 73.01] [color={rgb, 255:red, 0; green, 0; blue, 0 }  ][fill={rgb, 255:red, 0; green, 0; blue, 0 }  ][line width=0.75]      (0, 0) circle [x radius= 1.34, y radius= 1.34]   ;
\draw [shift={(271.63,173.2)}, rotate = 73.01] [color={rgb, 255:red, 0; green, 0; blue, 0 }  ][fill={rgb, 255:red, 0; green, 0; blue, 0 }  ][line width=0.75]      (0, 0) circle [x radius= 1.34, y radius= 1.34]   ;
\draw  [dash pattern={on 0.84pt off 2.51pt}]  (111.68,189.9) -- (118.72,175.37) ;
\draw [shift={(118.72,175.37)}, rotate = 295.82] [color={rgb, 255:red, 0; green, 0; blue, 0 }  ][fill={rgb, 255:red, 0; green, 0; blue, 0 }  ][line width=0.75]      (0, 0) circle [x radius= 1.34, y radius= 1.34]   ;
\draw [shift={(111.68,189.9)}, rotate = 295.82] [color={rgb, 255:red, 0; green, 0; blue, 0 }  ][fill={rgb, 255:red, 0; green, 0; blue, 0 }  ][line width=0.75]      (0, 0) circle [x radius= 1.34, y radius= 1.34]   ;
\draw  [dash pattern={on 0.84pt off 2.51pt}]  (479.63,154.2) -- (497.63,160.2) ;
\draw [shift={(497.63,160.2)}, rotate = 18.43] [color={rgb, 255:red, 0; green, 0; blue, 0 }  ][fill={rgb, 255:red, 0; green, 0; blue, 0 }  ][line width=0.75]      (0, 0) circle [x radius= 1.34, y radius= 1.34]   ;
\draw [shift={(479.63,154.2)}, rotate = 18.43] [color={rgb, 255:red, 0; green, 0; blue, 0 }  ][fill={rgb, 255:red, 0; green, 0; blue, 0 }  ][line width=0.75]      (0, 0) circle [x radius= 1.34, y radius= 1.34]   ;
\draw  [dash pattern={on 0.84pt off 2.51pt}]  (512.63,173.2) -- (517.67,189.67) ;
\draw [shift={(517.67,189.67)}, rotate = 73.01] [color={rgb, 255:red, 0; green, 0; blue, 0 }  ][fill={rgb, 255:red, 0; green, 0; blue, 0 }  ][line width=0.75]      (0, 0) circle [x radius= 1.34, y radius= 1.34]   ;
\draw [shift={(512.63,173.2)}, rotate = 73.01] [color={rgb, 255:red, 0; green, 0; blue, 0 }  ][fill={rgb, 255:red, 0; green, 0; blue, 0 }  ][line width=0.75]      (0, 0) circle [x radius= 1.34, y radius= 1.34]   ;

\draw (162,104) node [anchor=north west][inner sep=0.75pt]  [font=\footnotesize] [align=left] {$\displaystyle n-k+1$};
\draw (226,26) node [anchor=north west][inner sep=0.75pt]  [font=\footnotesize] [align=left] {$\displaystyle k-1$};
\draw (186,129) node [anchor=north west][inner sep=0.75pt]  [font=\footnotesize] [align=left] {$\displaystyle T^{*}$};
\draw (424,128) node [anchor=north west][inner sep=0.75pt]  [font=\footnotesize] [align=left] {$\displaystyle T_{1}^{*}$};
\draw (186,212) node [anchor=north west][inner sep=0.75pt]  [font=\footnotesize] [align=left] {$\displaystyle T_{2}^{*}$};
\draw (427,210) node [anchor=north west][inner sep=0.75pt]  [font=\footnotesize] [align=left] {$\displaystyle T_{3}^{*}$};
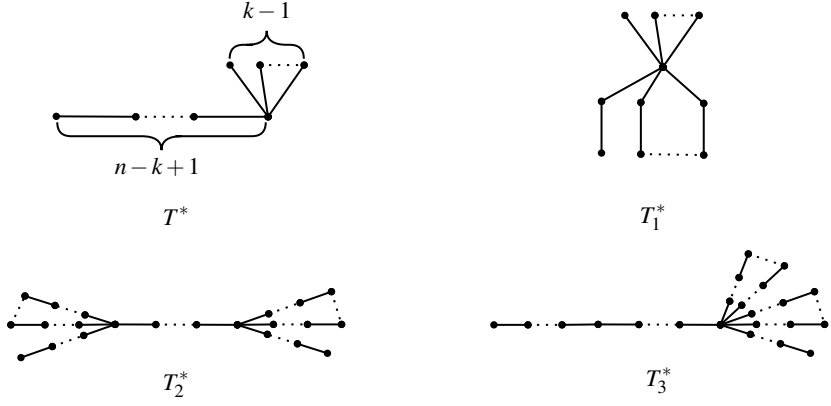
\draw (202,246) node [anchor=north west][inner sep=0.75pt]  [font=\small] [align=left] {Fig. 1. The graphs $\displaystyle T^{*}$, $\displaystyle T_{1}^{*}$, $\displaystyle T_{2}^{*}$ and $\displaystyle T_{3}^{*}$.};

\end{tikzpicture}
\end{center}

\section{Proof of Theorem \ref{th1}}

In this section, we will introduce some significant tools and then give the proof of Theorem \ref{th1}. Let $P_n$ be the path of order $n$, and $T_{a,b}$ denote the double star obtained by adding $a$ pendant vertices to one end vertex of $P_2$ and $b$ pendant vertices to the other. Throughout the discussion in \cite{D. Li}, we have the following remark.

\begin{remark}\label{rek5.11}
Let $\Sigma=(K_n, T^-_{s,t})$  be a signed complete graph and $n=s+t+2$. Then $\lambda_1(A((K_n, T^-_{s,t})))<\lambda_1(A((K_n, T^-_{s-1,t+1})))<\lambda_1(A((K_n, T^-_{1,n-3})))$ for $s\geq2$.
\end{remark}

\begin{lemma}\cite{T. Koledin} \label{lm2.01} Let $r,s,t$ and $u$ be distinct vertices of a signed graph $\Sigma$ and $X=(x_1,x_2,\ldots,x_n)^T$ be an eigenvector corresponding to $\lambda_1(A(\Sigma))$. Then
\begin{enumerate}
  \item [(i)] let $\Sigma'$ be obtained from $\Sigma$ by reversing the sign
of the positive edge $rs$ and the negative edge $rt$. If
$$\left\{
    \begin{array}{ll}
       x_r\geq0,x_t\geq x_s~~or\\
      x_r\leq0,x_t\leq x_s,
    \end{array}
  \right.$$
then $\lambda_1(A(\Sigma))\leq\lambda_1(A(\Sigma'))$. If at least one inequality for the entries of $X$ is strict, then $\lambda_1(A(\Sigma))<\lambda_1(A(\Sigma'))$;
  \item [(ii)] let $\Sigma'$ be obtained from $\Sigma$ by reversing the sign of
the positive edge $rs$ and the negative edge $tu$. If $x_rx_s\leq x_tx_u$, then $\lambda_1(A(\Sigma))\leq\lambda_1(A(\Sigma'))$. If at least one of the entries $x_r,x_s,x_t,x_u$ is distinct from zero, then $\lambda_1(A(\Sigma))<\lambda_1(A(\Sigma'))$.
\end{enumerate}
\end{lemma}

\textbf{Now we begin to prove Theorem \ref{th1}}. If $k=n-1$, then $T\cong K_{1,n-1}$. If $k=n-2$, then $T\cong T_{s,t}$, and $\lambda_1(A((K_n, T^-_{s,t})))\leq\lambda_1(A((K_n, T^-_{1,n-3})))$ with equality if and only if $T^-_{s,t}\cong T^-_{1,n-3}$ by Remark \ref{rek5.11}. If $k=2$, then $T\cong P_n$. Thus, we only need to consider $3\leq k\leq n-3$. Obviously, $(K_n,T^-)$ is unbalanced since the only balanced $(K_n,T^-)$ is the one with $T\cong K_{1,n-1}$.
And $3\leq\Delta(T)\le k$.
Let $X=(x_1,x_2,\ldots,x_n)^T$ be a unit eigenvector corresponding to $\lambda_1(A(\Sigma))$.
By arranging the vertices of $\Sigma$ appropriately, we can assume that $V(\Sigma)=\{v_1,v_2,\ldots,v_n\}$ such that $x_1\leq x_2\leq \cdots\leq x_n$, where $x_i$ corresponds to the vertex $v_i$ for $1\leq i\leq n$. Note that $-X$ must be an eigenvector of $\Sigma$ if $X$ is an eigenvector. Then we divide the proof into the following two cases.

\noindent{\bf{$\underline{\mbox{Case 1. }}$}} There exists a nonnegative eigenvector. Obviously, $x_n>0$ since $X\neq\mathbf{0}$.

\begin{claim}\label{claim21}
Up to replacing $v_1$ with another vertex $v_i$ such that $x_i=x_1$, we have $d_T(v_1)=\Delta(T)$.\end{claim}

It is trivial when $d_T(v_1)=\Delta(T)$. If $d_T(v_1)<\Delta(T)$, then there is a vertex $v_t\in V(T)$ such that $d_T(v_t)=\Delta(T)$. Let $P$ be the unique path in $T$ between $v_1$ and $v_t$. Then we assert that $x_t=x_1$. Otherwise, we will divide into the following two cases. If $d_T(v_1)=1$, then we can construct a new unbalanced signed graph $\Sigma'$ by reversing the signs of the positive edges $v_1w$ and the negative edges $v_tw$ whose negative edges also form a spanning tree with $k$ pendent vertices, where $w\in N_T(v_t)\setminus V(P)$. Thus,
\begin{align*}
\lambda_1(A(\Sigma'))-\lambda_1(A(\Sigma))&\geq X^T(A(\Sigma')-A(\Sigma))X\\
\ &=4\sum_{w\in N_T(v_t)\setminus V(P)}x_w(x_t-x_1)\\
\ &\geq0.
\end{align*}
If $\lambda_1(A(\Sigma'))=\lambda_1(A(\Sigma))$, then $X$ is also an eigenvector of $A(\Sigma')$ corresponding to $\lambda_1(A(\Sigma'))$. However,
$$\lambda_1(A(\Sigma'))x_w-\lambda_1(A(\Sigma))x_w= 2(x_t-x_1)>0$$
for any $w\in N_T(v_t)\setminus V(P)$, a contradiction.
If $d_T(v_1)\geq2$, then we can construct a new unbalanced signed graph $\Sigma'$ by reversing the sign of the positive edge $v_1w$ and the negative edge $v_tw$ for any $w\in N_T(v_t)\setminus V(P)$ whose negative edges also form a spanning tree with $k$ pendent vertices such that $\lambda_1(A(\Sigma))<\lambda_1(A(\Sigma'))$ by Lemma \ref{lm2.01}, a contradiction. Thus, $x_t=x_1$. And then we just need to exchange the subscripts $v_t$ and $v_1$, as desired.

Claim \ref{claim21} means that $3\leq d_T(v_1)=\Delta(T)\leq k$. Note that there are at least two vertices of $T$ with degree greater than three if $\Delta(T)\leq k-1$. Next, we will show that $d_T(v_1)=\Delta(T)=k$.

\begin{claim}\label{claim211} $\Delta(T)=k$.\end{claim}

Otherwise, there is a positive integer $t>1$ such that $d_T(v_t)\ge3$. Let $P$ be the path from $v_1$ to $v_t$ in $T$. Since $d_T(v_t)\ge 3$, there is a vertex $v_s$ such that $v_s \notin V(P)$ and $v_sv_t\in E(T)$. Then we first assert that $x_t=x_1$. Otherwise, we can construct a new unbalanced signed graph $\Sigma'$ by reversing the sign of the positive edge $v_1v_s$ and the negative edge $v_tv_s$ whose negative edges also form a spanning tree with $k$ pendent vertices such that $\lambda_1(A(\Sigma))<\lambda_1(A(\Sigma'))$ by Lemma \ref{lm2.01}, a contradiction. Let $N_1=N_T(v_1)\setminus V(P)$ and $N_t=N_T(v_t)\setminus V(P)$. Then $x_u=0$ for any $u\in N_1\cup N_t$ and $x_1=x_t=0$. Next, we will divide into the following three cases. If $v_1v_t\in E(T)$, then
$$0=\lambda_1(A(\Sigma))x_1=\sum_{u\in V(G)\setminus N_1\cup N_t\cup\{v_1,v_t\}}x_u,$$
which implies $x_u=0$ for any $u\in V(G)\setminus N_1\cup N_t\cup\{v_1,v_t\}$, and then $X=\mathbf{0}$, a contradiction. If $N_T(v_1)\cap N_T(v_t)=\{v_r\}$, i.e., $d_T(v_1,v_t)=2$, then we claim that $d_T(v_r)=2$. Otherwise, let $u\in N_r=N_T(v_r)\setminus \{v_1,v_t\}$. If $x_r>0$, then we can construct a new unbalanced signed graph $\Sigma'$ by reversing the sign of the positive edge $uv_1$ and the negative edge $uv_r$ whose negative edges also form a spanning tree with $k$ pendent vertices such that $\lambda_1(A(\Sigma))<\lambda_1(A(\Sigma'))$ by Lemma \ref{lm2.01}, a contradiction. So, $x_r=0$. However,
$$0=\lambda_1(A(\Sigma))x_1=\sum_{u\in V(G)\setminus (N_T[v_1]\cup N_T[v_t])}x_u,$$
which implies $x_u=0$ for any $u\in V(G)\setminus N_1\cup N_t\cup\{v_1,v_t\}$, and then $X=\mathbf{0}$, a contradiction.
Thus, by characteristic equation we get
$$\left\{
    \begin{array}{ll}
      \lambda_1(A(\Sigma))x_r=\sum_{u\in V(G)\setminus (N_T[v_1]\cup N_T[v_t])}x_u,\\
      \lambda_1(A(\Sigma))x_1=-x_r+\sum_{u\in V(G)\setminus (N_T[v_1]\cup N_T[v_t])}x_u.
    \end{array}
  \right.
$$
It shows that $(\lambda_1(A(\Sigma))-1)x_r=0$. And then $\sum_{u\in V(G)\setminus (N_T[v_1]\cup N_T[v_t])}x_u=x_r=0$ since $\lambda_1(A(\Sigma))>1$. So $X=\mathbf{0}$, a contradiction. Now we assume that $N_T(v_1)\cap N_T(v_t)=\emptyset$, i.e., $d_T(v_1,v_t)\geq3$. Let $P=v_1v_f\cdots v_hv_t$. By the similar discussion, we have $d_T(u)=2$ for any $u\in V(P)\setminus \{v_1,v_t\}$. By characteristic equation we get
$$\left\{
    \begin{array}{ll}
      \lambda_1(A(\Sigma))x_1=-x_f+x_h+\sum_{u\in V(G)\setminus (N_T[v_1]\cup N_T[v_t])}x_u,\\
      \lambda_1(A(\Sigma))x_t=x_f-x_h+\sum_{u\in V(G)\setminus (N_T[v_1]\cup N_T[v_t])}x_u.
    \end{array}
  \right.
$$
Thus, $x_f=x_h$ and $\sum_{u\in V(G)\setminus (N_T[v_1]\cup N_T[v_t])}x_u=0$. And then $x_f=x_h=0$. So $X=\mathbf{0}$, a contradiction.

By Claim \ref{claim211}, $d_T(v_i)\le 2$ for any positive integer $i\in [2,n]$.
Let $p$ and $q$ be the minimum and maximum subscripts such that $\{v_p,v_q\}\subset N_T(v_1)$, respectively. It is obvious that $x_p=\mbox{min}_{w\in N_T(v_1)}x_w$ and $x_q=\mbox{max}_{w\in N_T(v_1)}x_w$.

\begin{claim}\label{claim212} Up to replacing $v_p$ with another vertex $v_i$ such that $x_i=x_p$, we have $d_T(v_p)=2$.\end{claim}

It is trivial when $d_T(v_p)=2$. Now, we consider $d_T(v_p)=1$. Note that $\mid V(T)\setminus N_T[v_1]\mid\geq2$ since $k\leq n-3$, i.e., $n\geq k+1+2$. Let $v_s\in N_T(v_1)\setminus \{v_p\}$, $v_t\in V(T)\setminus N_T[v_1]$ and $v_sv_t\in E(T)$. Then we assert that $x_s=x_p$. Otherwise, we can construct a new unbalanced signed graph $\Sigma'$ by reversing the sign of the positive edge $v_pv_t$ and the negative edge $v_sv_t$ whose negative edges also form a spanning tree with $k$ pendent vertices such that $\lambda_1(A(\Sigma))<\lambda_1(A(\Sigma'))$ by Lemma \ref{lm2.01}, a contradiction. Thus, we just need to exchange the subscripts $v_s$ and $v_p$, as desired.

\begin{claim}\label{claim213} $x_p>0$.\end{claim}

Otherwise, $x_1=x_p=0$. Let $v_tv_p\in E(T)$ by Claim \ref{claim212}$, S_1=\sum_{v\in N_T(v_1)\setminus \{v_p\}}x_v$ and $S_2=\sum_{v\in V(T)}x_v-S_1-x_1-x_p-x_t$. By characteristic equation we get
$$\left\{
    \begin{array}{ll}
      0=\lambda_1(A(\Sigma))x_1=-S_1+S_2+x_t,\\
      0=\lambda_1(A(\Sigma))x_p=S_1+S_2-x_t.
    \end{array}
  \right.
$$
Then $S_2=0$ and $x_t=S_1$.  And $(\lambda_1(A(\Sigma))-1)x_t=0$ by $\lambda_1(A(\Sigma))x_t=S_1=x_t$. So, $x_t=S_1=0$ since $\lambda_1(A(\Sigma))>1$. Thus, $X=\mathbf{0}$, a contradiction.

If $d_T(v)=1$ for any vertex $v\in N_T(v_1)\setminus \{v_p\}$, then $T=T^*$. Next, assume that there is a positive integer $r\in (p,q]$ such that $v_r\in N_T(v_1)$ and $d_T(v_r)=2$. Let $P_1$ be the unique path in $T$ which contains $v_1$ and $v_p$, $P_2$ be the unique path in $T$ which contains $v_1$ and $v_r$. Suppose that $N_T(v_p)\setminus \{v_1\}=\{v_{p'}\}$ and $v_s\in V(P_2)$ such that $d_T(v_s)=1$. If $x_p> x_s$, then we can construct a new unbalanced signed graph $\Sigma'$ by reversing the sign of the positive edge $v_sv_{p'}$ and the negative edge $v_pv_{p'}$ whose negative edges also form a spanning tree with $k$ pendent vertices such that $\lambda_1(A(\Sigma))<\lambda_1(A(\Sigma'))$ by Lemma \ref{lm2.01}, a contradiction. Let $v_{s'}v_s\in E(P_2)$. If $x_p\leq x_s$, then we can construct a new unbalanced signed graph $\Sigma'$ by reversing the signs of the positive edges $v_1v_s$ and $v_{s'}v_p$ and the negative edges $v_1v_p$ and $v_{s'}v_s$ whose negative edges also form a spanning tree with $k$ pendent vertices, and
\begin{align*}
\lambda_1(A(\Sigma'))-\lambda_1(A(\Sigma))&\geq X^T(A(\Sigma')-A(\Sigma))X\\
\ &=4(x_1-x_{s'})(x_p-x_s)\\
\ &\geq0.
\end{align*}
If $\lambda_1(A(\Sigma'))=\lambda_1(A(\Sigma))$, then $X$ is also the unit eigenvector corresponding to $\lambda_1(A(\Sigma'))$ and
$$\left\{
    \begin{array}{ll}
      0=(\lambda_1(A(\Sigma'))-\lambda_1(A(\Sigma)))x_1=2(x_p-x_s),\\
      0=(\lambda_1(A(\Sigma'))-\lambda_1(A(\Sigma)))x_p=2(x_1-x_{s'}),
    \end{array}
  \right.
$$
that is, $x_p=x_s$ and $x_1=x_{s'}$. Let $N=\{v_1,v_p,v_{p'},v_s,v_{s'}\}$ and $M=\sum_{v\in V(T)\setminus N}x_v$. Since
$$\left\{
    \begin{array}{ll}
      \lambda_1(A(\Sigma))x_p=-x_1-x_{p'}+x_s+x_{s'}+M,\\
      \lambda_1(A(\Sigma))x_s=x_1+x_{p'}+x_p-x_{s'}+M,
    \end{array}
  \right.
$$
we have
$$\left\{
    \begin{array}{ll}
      x_1=x_{s'}=x_{p'}=0,~~~~~~~~~~~~~~~~~~~~~~~~~~~~~~~~~~~~~~~~~~~~~~~~~~~~~~~~~~~~~~~~~~~~~~~~~~(1)\\
      (\lambda_1(A(\Sigma))-1)x_p=(\lambda_1(A(\Sigma))-1)x_s=M.~~~~~~~~~~~~~~~~~~~~~~~~~~~~~~~~~~(2)
    \end{array}
  \right.
$$
We assert that $v_r\neq v_{s'}$. Otherwise, $x_1=x_p=x_r=x_{p'}=x_s=0$, which implies that $X=\mathbf{0}$, a contradiction. If $d_T(v_q)=1$, then
$$\lambda_1(A(\Sigma))x_q=-x_1+x_p+x_{p'}+x_s+x_{s'}+M-x_q=2x_p+(\lambda_1(A(\Sigma))-1)x_p-x_q$$
by (1) and (2), that is, $x_p=x_q>0$. If $d_T(v_q)=2$, let $N_T(v_q)\setminus \{v_1\}=v_{q'}$, then
$$(\lambda_1(A(\Sigma))+1)x_q=-x_1-x_{q'}+x_p+x_{p'}+x_s+x_{s'}+M-x_{q'}=-2x_{q'}+(\lambda_1(A(\Sigma))+1)x_p,$$
that is, $2x_{q'}=(\lambda_1(A(\Sigma))+1)(x_p-x_q)\leq0$. Thus, $x_{q'}=0$ and $x_p=x_q>0.$
So,
$$0=\lambda_1(A(\Sigma))x_1=-kx_p+x_{p'}+x_s+x_{s'}+M-(k-1)x_p=(\lambda_1(A(\Sigma))-1-2(k-1))x_p,$$
it means that $\lambda_1(A(\Sigma))=2k-1$. Next, we will consider that following two cases.

\noindent{\bf{$\underline{\mbox{Subase 1.1. }}$}} $v_rv_{s'}\in E(T)$, i.e., $P_2=v_1v_rv_{s'}v_s$.
Note that
$$0=\lambda_1(A(\Sigma))x_{s'}=x_1+x_p+x_{p'}-x_s-x_{r}+M-x_r=(\lambda_1(A(\Sigma))-1)x_p-2x_p.$$
We have $\lambda_1(A(\Sigma))=3$, and then $k=2$, a contradiction.

\noindent{\bf{$\underline{\mbox{Subase 1.2. }}$}} $v_rv_{s'}\notin E(T)$. Let $\{v_s,v_{s''}\}\subset N_T(v_{s'})\cap V(P_2)$, i.e., $P_2=v_1v_r\cdots v_{s''}v_{s'}v_s$. Since
$$0=\lambda_1(A(\Sigma))x_{s'}=-x_{s''}-x_s+x_1+x_p+x_{p'}+M-x_{s''}=-2x_{s''}+(\lambda_1(A(\Sigma))-1)x_p,$$
we have $x_{s''}=(k-1)x_p$. It leads to $x_v=0$ for any $v\in V(T)\setminus (N_T(v_1)\cup\{x_s,x_{s''}\})$. So,
$0=\lambda_1(A(\Sigma))x_{p'}=(2k-2)x_{p}$, which means that $k=1$, a contradiction.

\noindent{\bf{$\underline{\mbox{Case 2. }}$}} There are no nonnegative eigenvectors. Fixed an eigenvector $X=(x_1,x_2,\ldots,x_n)^T$, let $V_+=\{v_i|x_{v_i}\geq0\}$ and $V_-=\{v_i|x_{v_i}<0\}$. Clearly, $V_+\neq\emptyset$, $V_-\neq\emptyset$ and there must exist a vertex $v$ such that  $x_v>0$.
Let $\mid V_+\mid=a$, $\mid V_-\mid=b$ and $x_1\leq\cdots\leq x_b<0\leq x_{b+1}\leq\cdots\leq x_n$. Note that $a+b=n$. For convenience, set $u_i=v_{n-(i-1)}$ and $y_i=x_{n-(i-1)}$ for $1\leq i\leq a$, then $x_1\leq\cdots\leq x_b<0\leq y_a\leq\cdots\leq y_1$ and $y_1>0$. Obviously, $u_i\in V_+$ for $1\leq i\leq a$.

\noindent{\bf{$\underline{\mbox{Subcase 2.1. }}$}} $a=1$ or $b=1$. By symmetry, we just consider $a=1$.

\begin{claim}\label{claim221}
$d_T(u_1)=\Delta(T)$.
\end{claim}

It is trivial when $d_T(u_1)=\Delta(T)$. If $d_T(u_1)<\Delta(T)$, then there is a vertex $v_t\in V(T)$ such that $d_T(v_t)=\Delta(T)$. Let $P$ be the unique path in $T$ between $u_1$ and $v_t$. If $d_T(u_1)=1$, then we can construct a new unbalanced signed graph $\Sigma'$ by reversing the signs of the positive edges $u_1w$ and the negative edges $v_tw$ whose negative edges also form a spanning tree with $k$ pendent vertices, where $w\in N_T(v_t)\setminus V(P)$. Thus,
\begin{align*}
\lambda_1(A(\Sigma'))-\lambda_1(A(\Sigma))&\geq X^T(A(\Sigma')-A(\Sigma))X\\
\ &=4\sum_{w\in N_T(v_t)\setminus V(P)}x_w(x_t-y_1)\\
\ &\geq0.
\end{align*}
If $\lambda_1(A(\Sigma'))=\lambda_1(A(\Sigma))$, then $X$ is also an eigenvector of $A(\Sigma')$ corresponding to $\lambda_1(A(\Sigma'))$. However,
$$\lambda_1(A(\Sigma'))x_w-\lambda_1(A(\Sigma))x_w= 2(x_t-y_1)<0$$
for any $w\in N_T(v_t)\setminus V(P)$, a contradiction.
If $d_T(u_1)\geq2$, then we can construct a new unbalanced signed graph $\Sigma'$ by reversing the sign of the positive edge $u_1w$ and the negative edge $v_tw$ for any $w\in N_T(v_t)\setminus V(P)$ whose negative edges also form a spanning tree with $k$ pendent vertices such that $\lambda_1(A(\Sigma))<\lambda_1(A(\Sigma'))$ by Lemma \ref{lm2.01}, a contradiction.

Claim \ref{claim221} means that $3\leq d_T(u_1)=\Delta(T)\leq k$. Note that there are at least two vertices of $T$ with degree greater than three if $\Delta(T)\leq k-1$. Next, we will show that $d_T(u_1)=\Delta(T)=k$.

\begin{claim}\label{claim222} 
$\Delta(T)=k$.
\end{claim}

Otherwise, there is a vertex $v_t\in V_-$ such that $d_T(v_t)\ge3$. Let $P$ be the path from $u_1$ to $v_t$ in $T$. Since $d_T(v_t)\ge 3$, there is a vertex $v_s\in V_-\setminus V(P)$ such that $v_sv_t\in E(T)$. We can construct a new unbalanced signed graph $\Sigma'$ by reversing the sign of the positive edge $u_1v_s$ and the negative edge $v_tv_s$ whose negative edges also form a spanning tree with $k$ pendent vertices such that $\lambda_1(A(\Sigma))<\lambda_1(A(\Sigma'))$ by Lemma \ref{lm2.01}, a contradiction.

By Claim \ref{claim222}, $d_T(v_i)\le 2$ for any $v_i\in V_-$.
Let $q$ be the maximum subscripts such that $v_q\in N_T(u_1)$. Evidently, $x_q=\mbox{max}_{w\in N_T(u_1)}x_w$.

\begin{claim}\label{claim223} 
For any vertex $v_i$ with $x_i<x_q$, we have $v_i\in N_T(u_1)$.
\end{claim}

Otherwise, there is a positive integer $t$ ($t<q$) such that $x_t<x_q$ and $u_1v_t\notin E(T)$. Let $P$ be the path between $u_1$ and $v_t$ in $T$. If $v_q\notin V(P)$, then there must be a vertex $v_s$ such that $v_s\in N_P(v_t)$ and we can construct a new unbalanced signed graph $\Sigma'$ by reversing the signs of the positive edges $u_1v_t$, $v_sv_q$ and the negative edges $u_1v_q$, $v_tv_s$ whose negative edges also form a spanning tree with $k$ pendent vertices. Thus,
\begin{align*}
\lambda_1(A(\Sigma'))-\lambda_1(A(\Sigma))&\geq X^T(A(\Sigma')-A(\Sigma))X\\
\ &=4(y_1-x_s)(x_q-x_t)\\
\ &>0,
\end{align*}
a contradiction. If $v_q\in V(P)$ and $d_T(v_t)=1$, then $P=u_1v_q\cdots v_t$ and we can construct a new unbalanced signed graph $\Sigma'$ by reversing the sign of the positive edge $u_1v_t$ and the negative edge $u_1v_q$ whose negative edges also form a spanning tree with $k$ pendent vertices such that $\lambda_1(A(\Sigma))<\lambda_1(A(\Sigma'))$ by Lemma \ref{lm2.01}, a contradiction.
If $v_q\in V(P)$ and $d_T(v_t)=2$, then there must be a vertex $v_s\notin V(P)$ such that $v_sv_t\in E(T)$. We can construct a new unbalanced signed graph $\Sigma'$ by reversing the signs of the positive edges $u_1v_t$, $v_qv_s$ and the negative edges $u_1v_q$, $v_tv_s$ whose negative edges also form a spanning tree with $k$ pendent vertices, then
\begin{align*}
\lambda_1(A(\Sigma'))-\lambda_1(A(\Sigma))&\geq X^T(A(\Sigma')-A(\Sigma))X\\
\ &=4(y_1-x_s)(x_q-x_t)\\
\ &>0,
\end{align*}
a contradiction.

Let $p$ be the minimum subscript such that $v_p\notin N_T(u_1)$. If $p>q$, then the subscripts of vertices of $N_T(v_1)$ are consecutive. If $1\leq p<q$, then $x_p=x_q$ by Claim \ref{claim223}. Let $r$ ($p<r\leq q$) be the minimum subscript such that $v_r\in N_T(u_1)$. Then $x_r=x_p$ and we exchange the subscripts of $v_r$ and $v_p$.  Keep doing this procedure until the subscripts of vertices of $N_T(u_1)$ are consecutive. Without loss of generality, it is not restrictive to assume that $N_T(u_1)=\{v_1, v_2,\dots ,v_k\}$ by Claims \ref{claim221}-\ref{claim223}.

\begin{claim}\label{claim224} 
$d_T(v_k)=2.$
\end{claim}

Otherwise, $d_T(v_k)=1.$ Recall that $3\leq\Delta(T)\le k$, thus there must exist vertices $v_i$ and $v_j$ such that $v_iv_j\in E(T)$, where $1\leq i\leq k-1$ and $j\geq k+1$. Note that $x_i\leq x_k\leq x_j<0$. Then we can construct a new unbalanced signed graph $\Sigma'$ by reversing the sign of the positive edge $v_kv_j$ and the negative edge $v_iv_j$ whose negative edges also form a spanning tree with $k$ pendent vertices such that $\lambda_1(A(\Sigma))<\lambda_1(A(\Sigma'))$ by Lemma \ref{lm2.01}, a contradiction.

By Claim \ref{claim224}, let $v_s\in N_T(v_k)\setminus \{u_1\}$. Then the following claim holds.

\begin{claim}\label{claim225} 
$x_s=x_{n-1}$.
\end{claim}

Otherwise, $x_s<x_{n-1}$. Let $P$ be the path between $v_k$ and $v_{n-1}$ in $T$. If $v_s\notin V(P)$, then there is a vertex $v_t\in N_T(v_{n-1})$ such that $v_tv_{n-1}\in E(P)$. Note that $x_t\leq x_k\leq x_s<x_{n-1}<0$. If $u_1v_t\in E(T)$, then we can construct a new unbalanced signed graph $\Sigma'$ by reversing the signs of the positive edges $v_kv_{n-1}$, $v_sv_t$ and the negative edges $v_kv_s$, $v_tv_{n-1}$ whose negative edges also form a spanning tree with $k$ pendent vertices and
\begin{align*}
\lambda_1(A(\Sigma'))-\lambda_1(A(\Sigma))&\geq X^T(A(\Sigma')-A(\Sigma))X\\
\ &=4[x_s(x_k-x_{n-1})+x_{n-1}(x_t-x_k)]\\
\ &>0,
\end{align*}
a contradiction. If $u_1v_t\notin E(T)$, then we can construct a new unbalanced signed graph $\Sigma'$ by reversing the sign of the positive edge $v_kv_{n-1}$ and the negative edge $v_kv_s$ whose negative edges also form a spanning tree with $k$ pendent vertices and
\begin{align*}
\lambda_1(A(\Sigma'))-\lambda_1(A(\Sigma))&\geq X^T(A(\Sigma')-A(\Sigma))X\\
\ &=4(x_k-x_t)(x_s-x_{n-1})\\
\ &\geq0.
\end{align*}
If $\lambda_1(A(\Sigma'))=\lambda_1(A(\Sigma))$, then $X$ is also an eigenvector of $A(\Sigma')$ corresponding to $\lambda_1(A(\Sigma'))$. However,
$$\lambda_1(A(\Sigma'))x_k-\lambda_1(A(\Sigma))x_k= 2(x_s-x_{n-1})<0,$$
a contradiction. Next, we assume that $v_s\in V(P)$. Note that $x_k\leq x_s<x_{n-1}<0$. If $d_T(v_{n-1})=1$, then we can construct a new unbalanced signed graph $\Sigma'$ by reversing the sign of the positive edge $v_kv_{n-1}$ and the negative edge $v_kv_s$ whose negative edges also form a spanning tree with $k$ pendent vertices such that $\lambda_1(A(\Sigma))<\lambda_1(A(\Sigma'))$ by Lemma \ref{lm2.01}, a contradiction. If $d_T(v_{n-1})=2$, let $v_t\in N_T(v_{n-1})$ such that $v_tv_{n-1}\in E(T)\setminus E(P)$. Thus, $x_k\leq x_t$ and we can construct a new unbalanced signed graph $\Sigma'$ by reversing the signs of the positive edges $v_kv_{n-1}$, $v_sv_t$ and the negative edges $v_kv_s$, $v_tv_{n-1}$ whose negative edges also form a spanning tree with $k$ pendent vertices and
\begin{align*}
\lambda_1(A(\Sigma'))-\lambda_1(A(\Sigma))&\geq X^T(A(\Sigma')-A(\Sigma))X\\
\ &=4(x_k-x_t)(x_s-x_{n-1})\\
\ &\geq0.
\end{align*}
If $\lambda_1(A(\Sigma'))=\lambda_1(A(\Sigma))$, then $X$ is also an eigenvector of $A(\Sigma')$ corresponding to $\lambda_1(A(\Sigma'))$. However,
$$\lambda_1(A(\Sigma'))x_k-\lambda_1(A(\Sigma))x_k= 2(x_s-x_{n-1})<0,$$
a contradiction.

If $s=n-1$, then $v_kv_{n-1}\in E(T)$. If $s<n-1$, then we just need to exchange the subscripts of $v_s$ and $v_{n-1}$ by Claim \ref{claim225}.
Without loss of generality, it is not restrictive to assume that $N_T(v_k)=\{u_1,v_{n-1}\}$.

\begin{claim}\label{claim226} 
$d_T(v_i)=1$ for $1\leq i\leq k-1$.
\end{claim}

Suppose by way of contradiction that there exist positive integers $t\in [1,k-1]$ and $s\in[k+1,n-2]$ such that $v_tv_s\in E(T)$. Let $P$ be the path which contains $u_1$ and $v_k$ in $T$, and $v_r\in V(P)$ such that $d_T(v_r)=1$. Obviously, $x_t\leq x_r$. Then we can construct a new unbalanced signed graph $\Sigma'$ by reversing the sign of the positive edge $v_sv_r$ and the negative edge $v_sv_t$ whose negative edges also form a spanning tree with $k$ pendent vertices such that $\lambda_1(A(\Sigma))<\lambda_1(A(\Sigma'))$ by Lemma \ref{lm2.01}, a contradiction.

By Claims \ref{claim221}-\ref{claim226}, we obtain that $T=T^*$ immediately.

\noindent{\bf{$\underline{\mbox{Subcase 2.2. }}$}} $a\geq2$ and $b\geq2$. Recall that $-X$ must be an eigenvector of $\Sigma$ if $X$ is an eigenvector. So, without loss of generality, assume that $u_r\in V_+$ is the vertex such that $d_T(u_r)=\Delta(T)$. Let $N_r^+=N_T(u_r)\cap V_+$ and
$N_r^-=N_T(u_r)\cap V_-$.
\begin{claim}\label{claim2221}
If $N_r^+\neq \emptyset$, then $N_r^-\neq \emptyset$.
\end{claim}

Assume to the contrary that $N_r^-=\emptyset$. Then there is a vertex $u_t\in N_r^+$ such that $P=u_ru_t\cdots v_1$ is the unique path between $u_r$ and $v_1$ in $T$. Since $3\leq\Delta(T)\le k$, $\mid N_r^+\setminus \{u_t\}\mid\neq \emptyset$. If $d_T(v_1)=1$, then we can construct a new unbalanced signed graph $\Sigma'$ by reversing the sign of the positive edge $v_1w$ and the negative edge $u_rw$ for all $w\in N_r^+\setminus \{u_t\}$ whose negative edges also form a spanning tree with $k$ pendent vertices and
\begin{align*}
\lambda_1(A(\Sigma'))-\lambda_1(A(\Sigma))&\geq X^T(A(\Sigma')-A(\Sigma))X\\
\ &=4\sum_{w\in N_r^+\setminus \{u_t\}}y_w(y_r-x_1)\\
\ &\geq0.
\end{align*}
If $\lambda_1(A(\Sigma'))=\lambda_1(A(\Sigma))$, then $X$ is also an eigenvector of $A(\Sigma')$ corresponding to $\lambda_1(A(\Sigma'))$. However, for any $w\in N_r^+\setminus \{u_t\}$ we have
$$\lambda_1(A(\Sigma'))y_w-\lambda_1(A(\Sigma))y_w= 2(y_r-x_1)>0,$$
a contradiction.
If $d_T(v_1)\geq 2$, let $w\in N_r^+\setminus \{u_t\}$, then we can construct a new unbalanced signed graph $\Sigma'$ by reversing the sign of the positive edge $v_1w$ and the negative edge $u_rw$ whose negative edges also form a spanning tree with $k$ pendent vertices such that $\lambda_1(A(\Sigma))<\lambda_1(A(\Sigma'))$ by Lemma \ref{lm2.01}, a contradiction.

By Claim \ref{claim2221}, we divide the proof into the following two cases.

\noindent{\bf{$\underline{\mbox{Subcase 2.2.1 }}$}} $N_r^+\neq \emptyset$. Note that $N_r^-\neq \emptyset$ by Claim \ref{claim2221}.

\begin{claim}\label{claim2222}
$d_T(u)=1$ for any $u\in N_r^-$.
\end{claim}

Otherwise, suppose that $v_i\in N_r^-$ such that $d_T(v_i)\geq 2$. We assert that $\mid N_r^-\mid\geq 2$. Otherwise, let $w\in N_r^+$, then we can construct a new unbalanced signed graph $\Sigma'$ by reversing the sign of the positive edge $v_iw$ and the negative edge $u_rw$ whose negative edges also form a spanning tree with $k$ pendent vertices such that $\lambda_1(A(\Sigma))<\lambda_1(A(\Sigma'))$ by Lemma \ref{lm2.01}, a contradiction.
Now, we can construct a new unbalanced signed graph $\Sigma'$ by reversing the sign of the positive edge $v_iw$ and the negative edge $u_rw$ for all $w\in N_r^+$ whose negative edges also form a spanning tree with $k$ pendent vertices and
\begin{align*}
\lambda_1(A(\Sigma'))-\lambda_1(A(\Sigma))&\geq X^T(A(\Sigma')-A(\Sigma))X\\
\ &=4\sum_{w\in N_r^+}y_w(y_r-x_i)\\
\ &\geq0.
\end{align*}
If $\lambda_1(A(\Sigma'))=\lambda_1(A(\Sigma))$, then $X$ is also an eigenvector of $A(\Sigma')$ corresponding to $\lambda_1(A(\Sigma'))$. However, for any $w\in N_r^+$ we have
$$\lambda_1(A(\Sigma'))y_w-\lambda_1(A(\Sigma))y_w= 2(y_r-x_i)>0,$$
a contradiction.

\begin{claim}\label{claim2223}
$y_r\geq y_i$ for any $u_i\in N_r^+$.
\end{claim}

Assume to the contrary that there is a vertex $u_s\in N_r^+$ such that $y_r<y_s$. If $d_T(u_s)=1$, then by characteristic equation we get
$$\left\{
    \begin{array}{ll}
      (\lambda_1+1)(y_r+y_s)=2\sum\limits_{w\in (V(T)\setminus N_T[u_r])}x_w,\\
      (\lambda_1-1)(y_s-y_r)=2\sum\limits_{w\in (N_T(u_r)\setminus \{u_s\})}x_w.
    \end{array}
  \right.
$$
These mean that $\sum\limits_{w\in (V(T)\setminus N_T[u_r])}x_w\geq0$ and $\sum\limits_{w\in (N_T(u_r)\setminus \{u_s\})}x_w>0$.
By Claim \ref{claim2222}, for any vertex $v_i\in N_r^-$ we have
$$(\lambda_1+1)x_i+y_r-y_s-\sum\limits_{w\in (N_T(u_r)\setminus \{u_s\})}x_w=\sum\limits_{w\in (V(T)\setminus N_T[u_r])}x_w,$$
which means that $\sum\limits_{w\in (V(T)\setminus N_T[u_r])}x_w<0$, a contradiction. Now, we assume that $d_T(u_s)\geq2$. Recall that $N_r^-\neq \emptyset$ and let $u\in N_r^-$. Then we can construct a new unbalanced signed graph $\Sigma'$ by reversing the sign of the positive edge $u_su$ and the negative edge $u_ru$ whose negative edges also form a spanning tree with $k$ pendent vertices such that $\lambda_1(A(\Sigma))<\lambda_1(A(\Sigma'))$ by Lemma \ref{lm2.01}, a contradiction.

\begin{claim}\label{claim2224}
$y_r=y_1$.
\end{claim}
Otherwise, $y_r<y_1$. We first assume that $d_T(u_1)=1$ and $u_1u_r\in E(T)$. Then there is a vertex $u_s\in N_r^+$ such that $d_T(u_s)\geq2$ since $T\neq K_{1,n-1}$. And we can construct a new unbalanced signed graph $\Sigma'$ by reversing the sign of the positive edge $u_sw$ and the negative edge $u_rw$ for all $w\in N_r^+\setminus \{u_s\}$ and reversing the sign of the positive edge $u_1v$ and the negative edge $u_rv$ for all $v\in N_r^-$ whose negative edges also form a spanning tree with $k$ pendent vertices and
\begin{align*}
\lambda_1(A(\Sigma'))-\lambda_1(A(\Sigma))&\geq X^T(A(\Sigma')-A(\Sigma))X\\
\ &=4\sum_{w\in N_r^+\setminus \{u_s\}}y_w(y_r-y_s)+4\sum_{v\in N_r^-}x_v(y_r-y_1)\\
\ &>0,
\end{align*}
a contradiction.
Now we assume that $d_T(u_1)=1$ and $u_1u_r\notin E(T)$. Let $P=u_1\cdots u_tu_r$ be the path between $u_1$ and $u_r$ in $T$. Then we can construct a new unbalanced signed graph $\Sigma'$ by reversing the sign of the positive edge $u_tw$ and the negative edge $u_rw$ for all $w\in N_r^+\setminus \{u_t\}$ and reversing the sign of the positive edge $u_1v$ and the negative edge $u_rv$ for all $v\in N_r^-$ whose negative edges also form a spanning tree with $k$ pendent vertices and
\begin{align*}
\lambda_1(A(\Sigma'))-\lambda_1(A(\Sigma))&\geq X^T(A(\Sigma')-A(\Sigma))X\\
\ &=4\sum_{w\in N_r^+\setminus \{u_t\}}y_w(y_r-y_s)+4\sum_{v\in N_r^-}x_v(y_r-y_1)\\
\ &>0,
\end{align*}
a contradiction. Finally, let $d_T(u_1)\geq2$ and $v\in N_r^-$. Note that $d_T(u_r)=\Delta(T)\geq3$. Then we can construct a new unbalanced signed graph $\Sigma'$ by reversing the sign of the positive edge $u_1v$ and the negative edge $u_rv$ whose negative edges also form a spanning tree with $k$ pendent vertices such that $\lambda_1(A(\Sigma))<\lambda_1(A(\Sigma'))$ by Lemma \ref{lm2.01}, a contradiction.

It is not restrictive to assume that $d_T(u_1)=\Delta(T)\geq3$ by Claim \ref{claim2224}. Let $x_p=\mbox{min}_{w\in N_1^-}x_w$.

\begin{claim}\label{claim2225}
$x_p=x_1$.
\end{claim}
If $p=1$, then the result holds obviously. Now we assume that $p\geq2$. Assume to the contrary that $x_p>x_1$. By Claim \ref{claim2222}, there is a vertex $u_t\in N_r^+$ such that $P=v_1\cdots u_tu_1v_p$ is the path between $v_1$ and $v_p$ in $T$. If $v_1u_t\in E(T)$, then we can construct a new unbalanced signed graph $\Sigma'$ by reversing the signs of the positive edges $u_1v_1$, $u_tv_p$ and the negative edges $u_1v_p$, $u_tv_1$ whose negative edges also form a spanning tree with $k$ pendent vertices and
\begin{align*}
\lambda_1(A(\Sigma'))-\lambda_1(A(\Sigma))&\geq X^T(A(\Sigma')-A(\Sigma))X\\
\ &=4(y_1-y_t)(x_p-x_1)\\
\ &\geq0.
\end{align*}
If $\lambda_1(A(\Sigma'))=\lambda_1(A(\Sigma))$, then $X$ is also an eigenvector of $A(\Sigma')$ corresponding to $\lambda_1(A(\Sigma'))$. However,
$$\lambda_1(A(\Sigma'))y_1-\lambda_1(A(\Sigma))y_1= 2(x_p-x_1)>0,$$
a contradiction. If $v_1u_t\notin E(T)$, then there exists a vertex $w\in N_T(v_1)\cap E(P)$. If $w\in V_+$, then we can construct a new unbalanced signed graph $\Sigma'$ by reversing the signs of the positive edges $u_1v_1$, $wv_p$ and the negative edges $u_1v_p$, $wv_1$ whose negative edges also form a spanning tree with $k$ pendent vertices and
\begin{align*}
\lambda_1(A(\Sigma'))-\lambda_1(A(\Sigma))&\geq X^T(A(\Sigma')-A(\Sigma))X\\
\ &=4(y_1-y_w)(x_p-x_1)\\
\ &\geq0.
\end{align*}
If $\lambda_1(A(\Sigma'))=\lambda_1(A(\Sigma))$, then $X$ is also an eigenvector of $A(\Sigma')$ corresponding to $\lambda_1(A(\Sigma'))$. However,
$$\lambda_1(A(\Sigma'))y_1-\lambda_1(A(\Sigma))y_1= 2(x_p-x_1)>0,$$
a contradiction. If $w\in V_-$, then we can construct a new unbalanced signed graph $\Sigma'$ by reversing the signs of the positive edges $u_1v_1$, $wv_p$ and the negative edges $u_1v_p$, $wv_1$ whose negative edges also form a spanning tree with $k$ pendent vertices and
\begin{align*}
\lambda_1(A(\Sigma'))-\lambda_1(A(\Sigma))&\geq X^T(A(\Sigma')-A(\Sigma))X\\
\ &=4(y_1-x_w)(x_p-x_1)\\
\ &>0,
\end{align*}
a contradiction.

\begin{claim}\label{claim2226}
$d_T(u_1)=\Delta(T)=k$.
\end{claim}
Otherwise, $d_T(u_1)=\Delta(T)<k$. Then there is a vertex $w\in V\setminus \{u_1\}$ such that $d_T(w)\ge3$. Let $P$ be the path from $u_1$ to $w$ in $T$.
Assume that $w\in V_+$ at first. We assert that $N_T(w)\cap V_-=\emptyset$. Otherwise, let $v\in N_T(w)\cap V_-$. Then we can construct a new unbalanced signed graph $\Sigma'$ by reversing the sign of the positive edge $u_1v$ and the negative edge $wv$ whose negative edges also form a spanning tree with $k$ pendent vertices such that $\lambda_1(A(\Sigma))<\lambda_1(A(\Sigma'))$ by Lemma \ref{lm2.01}, a contradiction.
Since $d_T(w)\ge 3$, $N_T(w)\setminus V(P)\neq \emptyset$. Then we can construct a new unbalanced signed graph $\Sigma'$ by reversing the signs of the positive edge $uv_1$ and the negative edges $uw$ for all $u\in N_T(w)\setminus V(P)$ whose negative edges also form a spanning tree with $k$ pendent vertices and
\begin{align*}
\lambda_1(A(\Sigma'))-\lambda_1(A(\Sigma))&\geq X^T(A(\Sigma')-A(\Sigma))X\\
\ &=4(y_w-x_1)\sum_{u\in N_T(w)\setminus V(P)}y_u\\
\ &\geq0.
\end{align*}
If $\lambda_1(A(\Sigma'))=\lambda_1(A(\Sigma))$, then $X$ is also an eigenvector of $A(\Sigma')$ corresponding to $\lambda_1(A(\Sigma'))$. However,
$$\lambda_1(A(\Sigma'))y_u-\lambda_1(A(\Sigma))y_u= 2(y_w-x_1)>0$$
for any $u\in N_T(w)\setminus V(P)$, a contradiction. Now we assume that $w\in V_-$. Note that at least one of $(N_T(w)\setminus V(P))\cap V_+\neq\emptyset$ and $(N_T(w)\setminus V(P))\cap V_-\neq\emptyset$ holds since $d_T(w)\ge3$. If $(N_T(w)\setminus V(P))\cap V_-\neq\emptyset$, then
we can construct a new unbalanced signed graph $\Sigma'$ by reversing the signs of the positive edges $vv_1$, $uu_1$ and the negative edges $vw$, $uw$ for all $u\in (N_T(w)\setminus V(P))\cap V_-$ and $v\in (N_T(w)\setminus V(P))\cap V_+$ whose negative edges also form a spanning tree with $k$ pendent vertices and
\begin{align*}
\lambda_1(A(\Sigma'))-\lambda_1(A(\Sigma))&\geq X^T(A(\Sigma')-A(\Sigma))X\\
\ &=4(x_w-x_1)\sum_{v\in (N_T(w)\setminus V(P))\cap V_+}y_v+4(x_w-y_1)\sum_{u\in (N_T(w)\setminus V(P))\cap V_-}x_u\\
\ &>0,
\end{align*}
a contradiction.
If $(N_T(w)\setminus V(P))\cap V_-=\emptyset$, then
we can construct a new unbalanced signed graph $\Sigma'$ by reversing the signs of the positive edges $u_tv_1$, $vv_1$ and the negative edges $u_tu_1$, $vw$ for all $v\in (N_T(w)\setminus V(P))\cap V_+$ whose negative edges also form a spanning tree with $k$ pendent vertices and
\begin{align*}
\lambda_1(A(\Sigma'))-\lambda_1(A(\Sigma))&\geq X^T(A(\Sigma')-A(\Sigma))X\\
\ &=4[y_t(y_1-x_1)+(x_w-x_1)\sum_{v\in (N_T(w)\setminus V(P))\cap V_+}y_v]\\
\ &\geq0.
\end{align*}
If $\lambda_1(A(\Sigma'))=\lambda_1(A(\Sigma))$, then $X$ is also an eigenvector of $A(\Sigma')$ corresponding to $\lambda_1(A(\Sigma'))$. However,
$$\lambda_1(A(\Sigma'))y_t-\lambda_1(A(\Sigma))y_t= 2(y_1-x_1)>0,$$
a contradiction.

Claim \ref{claim2226} shows that $d_T(u)\leq2$ for any $u\in V\setminus \{u_1\}$, which means that $T$ must be a starlike tree. Since $T\neq K_{1,n-1}$, there is a vertex $u_s\in N_T(u_1)\cap V_+$ such that $d_T(u_s)=2$. Let $w\in N_T(u_s)\setminus \{u_1\}$. We assert that $w\in V_+$. Otherwise, we can construct a new unbalanced signed graph $\Sigma'$ by reversing the signs of the positive edges $u_1w$, $v_1u_s$ and the negative edges $u_1u_s$, $u_sw$ whose negative edges also form a spanning tree with $k$ pendent vertices and
\begin{align*}
\lambda_1(A(\Sigma'))-\lambda_1(A(\Sigma))&\geq X^T(A(\Sigma')-A(\Sigma))X\\
\ &=4[y_1(y_s-x_w)+y_s(x_w-x_1)]\\
\ &>0,
\end{align*}
a contradiction. If $d_T(w)=2$, set $v\in N_T(w)\setminus \{u_s\}$, then either $v\in V_-$ or $v\in V_+$. If $v\in V_-$, then we can construct a new unbalanced signed graph $\Sigma'$ by reversing the signs of the positive edges $u_1w$, $v_1u_s$, $u_sv$ and the negative edges $u_1u_s$, $u_sw$, $wv$ whose negative edges also form a spanning tree with $k$ pendent vertices and
\begin{align*}
\lambda_1(A(\Sigma'))-\lambda_1(A(\Sigma))&\geq X^T(A(\Sigma')-A(\Sigma))X\\
\ &=4[(y_1-x_v)(y_s-x_w)+y_s(x_w-x_1)]\\
\ &>0,
\end{align*}
a contradiction. If $v\in V_+$, then we can construct a new unbalanced signed graph $\Sigma'$ by reversing the signs of the positive edges $u_1w$, $v_1u_s$, $u_sv$ and the negative edges $u_1u_s$, $u_sw$, $wv$ whose negative edges also form a spanning tree with $k$ pendent vertices and
\begin{align*}
\lambda_1(A(\Sigma'))-\lambda_1(A(\Sigma))&\geq X^T(A(\Sigma')-A(\Sigma))X\\
\ &=4[(y_1-y_v)(y_s-x_w)+y_s(x_w-x_1)]\\
\ &\geq0.
\end{align*}
If $\lambda_1(A(\Sigma'))=\lambda_1(A(\Sigma))$, then $X$ is also an eigenvector of $A(\Sigma')$ corresponding to $\lambda_1(A(\Sigma'))$. However,
$$\lambda_1(A(\Sigma'))y_1-\lambda_1(A(\Sigma))y_1= 2(y_s-x_w)>0,$$ a contradiction. These show that $d_T(w)=1$ and $T\cong T^*_1$. Let $\{u_s,u_t,u_p,u_q\}\subset V_+$ such that $\{u_1u_s,u_1u_t,u_su_p,u_tu_q\}\subset E(T)$, then we can construct a new unbalanced signed graph $\Sigma'$ by reversing the sign of the positive edge $u_pu_t$ and the negative edge $u_1u_t$ whose negative edges also form a spanning tree with $k$ pendent vertices such that $\lambda_1(A(\Sigma))<\lambda_1(A(\Sigma'))$ by Lemma \ref{lm2.01}. Continue the above procedure inductively, we get that $\lambda_1(A((\Sigma,T^*_1)))<\lambda_1(A((\Sigma',T^*)))$.

\noindent{\bf{$\underline{\mbox{Subcase 2.2.2 }}$}} $N_r^+= \emptyset$. Note that $\mid N_r^-\mid\geq3$ since $d_T(u_r)=\Delta(T)\geq3$.

\begin{claim}\label{claim2227}
$y_r=y_1$.
\end{claim}
Otherwise, $y_r<y_1$. Let $P=u_1\cdots v_su_r$ be the unique path between $u_1$ and $u_r$ in $T$, where $v_s\in N_r^-$. If $d_T(u_1)=1$, then we can construct a new unbalanced signed graph $\Sigma'$ by reversing the sign of the positive edge $u_1w$ and the negative edge $u_rw$ for all $w\in N_r^-\setminus \{v_s\}$ whose negative edges also form a spanning tree with $k$ pendent vertices and
\begin{align*}
\lambda_1(A(\Sigma'))-\lambda_1(A(\Sigma))&\geq X^T(A(\Sigma')-A(\Sigma))X\\
\ &=4\sum_{w\in N_r^-\setminus \{v_s\}}x_w(y_r-y_1)\\
\ &>0,
\end{align*}
a contradiction.
Now we assume that $d_T(u_1)\geq2$ and $w\in N_r^-\setminus \{v_s\}$. Then we can construct a new unbalanced signed graph $\Sigma'$ by reversing the sign of the positive edge $u_1w$ and the negative edge $u_rw$ whose negative edges also form a spanning tree with $k$ pendent vertices such that $\lambda_1(A(\Sigma))<\lambda_1(A(\Sigma'))$ by Lemma \ref{lm2.01}, a contradiction.

Let $x_p=\min_{w\in N_1^-}x_w$.

\begin{claim}\label{claim19}
$x_p=x_1$.
\end{claim}
If $p=1$, then the result holds obviously.
Now assume that $p\geq 2$. Assume on the contrary that $x_p>x_1$. There is a vertex $v_t\in N_1^-$ such that $P=v_1\cdots v_tu_1v_p$ or $P=v_1\cdots v_pu_1v_t$ is the path containing $v_1$, $v_p$ and $u_1$ in $T$. We first consider the case of $P=v_1\cdots v_tu_1v_p$. If $v_1v_t\in E(T)$, then we can construct a new unbalanced signed graph $\Sigma'$ by reversing the signs of the positive edges $u_1v_1, v_tv_p$ and the negative edges $u_1v_p, v_tv_1$ whose negative edges also form a spanning tree with $k$ pendent vertices and
\begin{align*}
\lambda_1(A(\Sigma'))-\lambda_1(A(\Sigma))&\geq X^T(A(\Sigma')-A(\Sigma))X\\
\ &=4(y_1-x_t)(x_p-x_1)\\
\ &>0,
\end{align*}
a contradiction.
If $v_1v_t\notin E(T)$, then there exists a vertex $w\in N_T(v_1)\cap E(P)$. If $w\in V_+$, then we can construct a new unbalanced signed graph $\Sigma'$ by reversing the signs of the positive edges $u_1v_1, wv_p$ and the negative edges $u_1v_p, wv_1$ whose negative edges also form a spanning tree with $k$ pendent vertices and
\begin{align*}
\lambda_1(A(\Sigma'))-\lambda_1(A(\Sigma))&\geq X^T(A(\Sigma')-A(\Sigma))X\\
\ &=4(y_1-y_w)(x_p-x_1)\\
\ &\geq 0.
\end{align*}
If $\lambda_1(A(\Sigma'))=\lambda_1(A(\Sigma))$, then $X$ is also an eigenvector of $A(\Sigma')$ corresponding to $\lambda_1(A(\Sigma'))$. However,
$$\lambda_1(A(\Sigma'))y_1-\lambda_1(A(\Sigma))y_1= 2(x_p-x_1)>0,$$ a contradiction. 
If $w\in V_-$,  then we can construct a new unbalanced signed graph $\Sigma'$ by reversing the signs of the positive edges $u_1v_1, wv_p$ and the negative edges $u_1v_p, wv_1$ whose negative edges also form a spanning tree with $k$ pendent vertices and
\begin{align*}
\lambda_1(A(\Sigma'))-\lambda_1(A(\Sigma))&\geq X^T(A(\Sigma')-A(\Sigma))X\\
\ &=4(y_1-x_w)(x_p-x_1)\\
\ &> 0,
\end{align*}
a contradiction.
Then we consider the case of $P=v_1\cdots v_pu_1v_t$. Similar to the above, we can get the contradiction.

It is not restrictive to assume that $d_T(u_1)=\Delta(T)\geq3$ by Claim \ref{claim2227}. 
Let $D_3=\{v\mid d_T(v)\geq3\}$. If $D_3\setminus\{u_1\}=\emptyset$, then $T$ is a starlike tree. Now we consider $D_3\setminus\{u_1\}\neq\emptyset$ and $w\in D_3\setminus\{u_1\}$. Let $P_1=u_1v_r\cdots w$ be the unique path between $u_1$ and $w$ in $T$.
If $w\in V_+$, then we assert that $x_v<0$ for any $v\in N_T(w)\setminus V(P_1)$. Assume to the contrary that $x_v\geq0$ and $v\in N_T(w)\setminus V(P_1)$, then we can construct a new unbalanced signed graph $\Sigma'$ by reversing the sign of the positive edge $vv_r$ and the negative edge $vw$ whose negative edges also form a spanning tree with $k$ pendent vertices such that $\lambda_1(A(\Sigma))<\lambda_1(A(\Sigma'))$ by Lemma \ref{lm2.01}, a contradiction. Similarly, we can obtain that $x_v>0$ for any $v\in N_T(w)\setminus V(P_1)$ if $w\in V_-$. And either $D_3\cap V_-=\emptyset$ or $D_3\cap (V_+\setminus \{u_1\})=\emptyset$. 
We first consider $D_3\cap V_-\neq\emptyset$ and $D_3\cap (V_+\setminus \{u_1\})=\emptyset$. 
Let $v_s\in D_3\cap V_-$. If $u_1v_s\notin E(T)$, let $P_2=u_1v_r\cdots vv_s$ be the unique path between $u_1$ and $v_s$ in $T$. Note that $v_r$ and $v$ may be the same vertex. Let $S= N_T(v_s)\setminus \{v\}$.
If $d_T(v_1)=2$, then for any $w\in S$, we can construct a new unbalanced signed graph $\Sigma'$ by reversing the signs of the positive edge $v_1w$ and the negative edge $v_sw$ whose negative edges also form a spanning tree with $k$ pendent vertices such that $\lambda_1(A(\Sigma))<\lambda_1(A(\Sigma'))$ by Lemma \ref{lm2.01}, a contradiction. 
If $d_T(v_1)=1$, then we can construct a new unbalanced signed graph $\Sigma'$ by reversing the signs of the positive edges $wv_1$ and the negative edges $wv_s$ for all $w\in S$ whose negative edges also form a spanning tree with $k$ pendent vertices and
\begin{align*}
\lambda_1(A(\Sigma'))-\lambda_1(A(\Sigma))&\geq X^T(A(\Sigma')-A(\Sigma))X\\
\ &=4\sum_{w\in S}y_w(x_s-x_1)\\
\ &\geq 0,
\end{align*}
If $\lambda_1(A(\Sigma'))=\lambda_1(A(\Sigma))$, then $X$ is also an eigenvector of $A(\Sigma')$ corresponding to $\lambda_1(A(\Sigma'))$. However,
$$\lambda_1(A(\Sigma'))x_1-\lambda_1(A(\Sigma))x_1= -2\Sigma_{w\in S}y_w<0,$$ a contradiction. 
If $u_1v_s\in E(T)$, we can similarly get the contradiction.
Thus, $D_3\cap V_-=\emptyset$. So, we only need to consider $D_3\cap (V_+\setminus \{u_1\})\neq\emptyset$. Similarly, we can get that $\mid D_3\cap (V_+\setminus \{u_1\})\mid=1$. Then $T\cong T^*_2$. 

Let $D_3\cap (V_+\setminus \{u_1\})=\{u_s\}$. Then $y_s\leq y_1$. For any $w\in N_T(u_s)$, we have $x_w<0$. Then we can construct a new unbalanced signed graph $\Sigma'$ by reversing the sign of the positive edge $u_1w$ and the negative edge $u_sw$ whose negative edges also form a spanning tree with $k$ pendent vertices such that $\lambda_1(A(\Sigma))<\lambda_1(A(\Sigma'))$ by Lemma \ref{lm2.01}, a contradiction. Continue the above procedure inductively, we get that $\lambda_1(A((\Sigma, T_2^*)))<\lambda_1(A((\Sigma', T_3^*)))$. One can easily observe that $\Delta(T_3^*)=k$. In this case, we have $\lambda_1(A((\Sigma', T_3^*)))\leq \lambda_1(A((\Sigma', T^*)))$. Therefore, $\lambda_1(A((\Sigma, T_2^*)))<\lambda_1(A((\Sigma', T^*)))$.
This completes the proof.

\section*{Acknowledgements}
The authors would like to show great gratitude to anonymous referees for their valuable suggestions which lead to an improvement of the original manuscript.

\end{document}